\documentclass[]{article}
\usepackage{bbm}
\usepackage{geometry}
\usepackage{amsbsy,amsfonts,amsmath,amssymb,enumerate,epsfig,graphicx,mathbbol,rotating,subfigure}

\newcommand{\1}{\textbf{1}}
\newtheorem{proposition}{Proposition}[section]
\newtheorem{lemma}{Lemma}[section]
\newtheorem{definition}{Definition}[section]

\newcommand{\Supp}{\mbox{Supp}}

\usepackage{bbm}

\newcommand{\sout}[1]{}

\newcommand{\xy}{
	\begin{pmatrix}x\\y\end{pmatrix}
}
 
\newcommand{\xystar}{
	\begin{pmatrix}x^*\\y^*\end{pmatrix}
}
 
\newcommand{\pileg}{ \tiny
	\begin{array}{c}
		y \geq x\\
		\xy \rightarrow \xystar
	\end{array}
}
\newcommand{\piled}{ \tiny
	\begin{array}{c}
		y < x\\
		\xy \rightarrow \xystar
	\end{array}
}

\begin{document}
\title{Density modification based reliability sensitivity analysis}

\author{P. Lema\^itre$^{\rm a}$ $^{\rm b}$ $^{\ast}$\thanks{$^\ast$Corresponding author. Email: lemaitre.paul.andre@gmail.com \vspace{6pt}}\and E. Sergienko$^{\rm c}$ $^{\rm d}$\and A. Arnaud$^{\rm a}$ \and N. Bousquet$^{\rm a}$ \and F. Gamboa$^{\rm d}$ \and B. Iooss$^{\rm a}$ $^{\rm d}$\\\vspace{6pt}
 \\
 $^{a}${\em{EDF R\&D, 6 Quai Watier - 78401 Chatou}};\\
 $^{b}${\em{INRIA Sud-Ouest, 351 cours de la lib\'eration - 33405 Talence}};\\
 $^{c}${\em{IFP EN, 1 avenue de Bois-Pr\'eau - 92852 Rueil-Malmaison}};\\
 $^{d}${\em{Institut de Math\'ematiques de Toulouse, 118 route de Narbonne - 31062 Toulouse}}}

\maketitle
\begin{abstract}
Sensitivity analysis of a numerical model, for instance simulating physical phenomena, is useful to quantify the influence of the inputs on the model responses.
This paper proposes a new sensitivity index, based upon the modification of the  probability density function (pdf) of the random inputs, when the quantity of interest is a failure probability (probability that a model
output exceeds a given threshold).
An input is considered influential if the input pdf modification leads to a broad change in the failure probability. 
These sensitivity indices can be computed using the sole set of simulations that has already been used to estimate the failure probability, thus limiting the number of calls to the numerical model.
In the case of a Monte Carlo sample, asymptotical properties of the indices are derived.
Based on Kullback-Leibler divergence, several types of input perturbations are introduced.
The relevance of this new sensitivity analysis method is analysed through three case studies.
\end{abstract}
%%%%%%%%%%%%%%%%%%%%%%%%%%%%%%%%%%%%%%%%%%%%%%%%%%%%%%%%%%%%%%%%%%%%%%%%%%%%%%%%%%%%%%%%%%%%%%%%%%%%%%%%%%%%%%%

%\begin{keywords}computer experiment; Kullback-Leibler; sensitivity analysis; structural reliability; uncertainty
%\end{keywords}

%%%%%%%%%%%%%%%%%%%%%%%%%%%%%%%%%%%%%%%%%%%%%%%%%%%%%%%%%%%%%%%%%%%%%%%%%%%%%%%%%%%%%%%%%%%%%%%%%%%%%%%%%%%%%%%

\section{Introduction}

In the context of structural reliability, computer models are used
in order to assess the safety of industrial systems relying on complex physical phenomena.
For instance, an electric operator would like to predict the level of a potential river 
flood in order to determine the height of a dyke preventing any disaster.
In this example, the computer model (simulating the hydraulic model) has some uncertain
input variables (flow rate, river length, water height, etc.), that are modelled
by random variables. In this paper, the computer code is a "black-box" deterministic
numerical model and the study focuses on one of its output. Due to the randomness of the model
inputs, this output is a random variable more or less sensitive to the uncertainty of the input variables.
 
Sensitivity analysis (SA) is a tool used to explore, understand and (partially) validate numerical models. It aims at explaining the outputs regarding the input
uncertainties (\cite{saltelli2002sensitivity}). We use the ``global SA'' definition
given by Saltelli \emph{et} al. \cite{saltelli2sensitivity} wherein
the whole variation range of the inputs is considered. The application
of such an approach can be model simplification (by removing irrelevant modelling
elements), input variables ranking or research prioritization.
There is a wide range of SA techniques, regarding what type of problem
the experimenter faces with (\cite{iooss2011}). For instance,
screening methods are to be applied when there is a large number of
inputs, and few models assumptions. From a quantitative point of view,
the most popular techniques are variance-based methods and the so-called Sobol' indices (\cite{saltelli2sensitivity,sobol1993sensitivity}).
These are based upon Hoeffding decomposition of $L^2$ function and functional variance decomposition \cite{Antoniadis}.

It should be noticed that most SA methods focus on real-valued continuous
numerical output variables. When the output is a binary value (e.g. when 
the numerical model returns ``faulty system'' or ``safe system''), SA 
techniques are underdeveloped.Some basic techniques can be quoted, such
as Monte-Carlo filtering (\cite{saltelli2sensitivity}) which consists in
 measuring differences between a ``safe'' sample and a ``faulty'' sample 
 via standard statistical tests.
In structural reliability analysis, some sensitivity factors resulting from the First or Second Order Reliability Methods (FORM/SORM, \cite{lemaire2009structural})
can also be used to classify the impact of the inputs on the failure
probability.
More recent works give methods combining the two objectives: estimating
a failure probability and assessing the influence of the input uncertainty
on this probability (\cite{morio2011influence,munoz2011adaptive}).

In this paper, a real-valued numerical model denoted by $G: \mathbb{R}^{d} \rightarrow \mathbb{R}$
is considered. This model may further be called the ``failure function''.
In practice, each run of $G$ can be CPU time consuming.
We are interested in the (rare) event
$G(\mathbf{X})<0$ (system failure) and in the complementary event $G(\mathbf{X})\geq0$
(system safe mode). $\mathbf{X}=(X_{1},...,X_{d})^{T}$
is a $d$-dimensional continuous random variable whose joint probability density
function (pdf) is denoted $f$. For $i=1,\cdots,d$, let  $f_{i}$ denotes the distribution of $X_i$ (the marginal pdf). 
We make the assumption that all components of $\mathbf{X}$
are independent. The quantity of interest is the system failure probability:
\begin{equation*}
P=\int \1_{\{G(\mathbf{x})<0\}}f(\mathbf{x}) d\mathbf{x}.
\end{equation*}
The aim of this work is the quantification of the influence of each variable $X_{i}$ on this probability. 

Let us ask the question: what are the engineer’s motivations when he perform a SA on his/her black-box model that produces a binary response? We provided an overview of the ”general objectives” of SA: variable ranking, model simplification, model understanding. But from our discussions with practitioners, we have identified three "engineer motivations":
	\begin{itemize}
	\item the practitioner wants to determine which are the inputs that impact the most the failure event – the inputs distributions being set and supposed to be perfectly known. This amounts to an absolute ranking objective.
	\item $P$ will be impacted by the choice of the input distributions; the engineer wants to assess the influence of this choice on the FP. Therefore the objective here is to quantify the sensitivity of the model output to the family or shape of the inputs.
	\item In practice, input distributions are estimated from data, thus leading to uncertainty on the values of the distribution parameters. The practitioner wants to assess the influence of the distribution parameters on $P$. Therefore the objective here is to test the sensitivity of the model to the parameters of the inputs.
	\end{itemize}

In most studies, sensitivity indices for failure probabilities are defined in strong correspondence with a given method
of estimation (e.g. \cite{lemaire2009structural,munoz2011adaptive}). Their interpretation is consequently limited.
We propose in this article to define new generic sensitivity indices.
Our sensitivity index is based upon density modification, and is adapted to failure probabilities.  A methodology to estimate such indices is derived.   
For simplicity reasons, a classical Monte Carlo framework is considered in the following. Additionally, the sensitivity index can be computed using the sole
set of simulations that has already been used to estimate the failure probability $P$, thus limiting the number of calls to the numerical model.

The outline of the article is the following: first we define a generic strategy of input perturbation in Section \ref{sec:Methodologies-of-input}, based upon maximum entropy rules.
We then present our index and its theoretical properties in Section \ref{sec:Definition}, altogether with the estimation methodology.  
The behaviour of the indices is examined in Section \ref{sec:Numerical-experiments} through numerical simulations in various complexity settings,
involving toy examples and a realistic case-study. Comparisons with two reference sensitivity analysis methods (FORM indices and Sobol' indices) 
highlight the relevance of the new indices in most situations. 
The main advantages and remaining issues are finally discussed in the last section of the article, that introduces avenues for future research. \\

%%%%%%%%%%%%%%%%%%%%%%%%%%%%%%%%%%%%%%%%%%%%%%%%%%%%%%%%%%%%%%%%%%%%%%%%%%%%%%%%%%%%%%%%%%%%%%%%%%%%%%%%%%%%%%%%%%%%%

\section{Methodologies of input perturbation\label{sec:Methodologies-of-input}}

Our sensitivity analysis method requires to define a perturbation for each input.
In general, and especially in preliminary reliability studies, there is no prior 
rule allowing to elicit a specialized perturbation for each input variable. 
We thus would like to propose a simple perturbation methodology, allowing the practitioner to answer the questions itemized in the Introduction.
Furthermore, we make the implicit hypothesis that the extreme values of the inputs lead to the (rare) failure event.

Given a unidimensional input variable $X_i$ with pdf $f_i$, let us call  $X_{i\delta}\sim f_{i\delta}$ the corresponding perturbed random input.
This perturbed input takes the place of the real random input $X_i$, in a sense of modelling error : what if the correct input were $X_{i\delta}$ instead of $X_i$?

More precisely, we suggest to define a perturbed input density $f_{i\delta}$ as the closest distribution to the original $f_{i}$ in the entropic sense and
under some constraints of perturbation. 
Information-theoretical arguments (\cite{cover2006elements}) led us to choose the Kullback-Leibler (KL) divergence  between $f_{i\delta}$ and $f_{i}$
 as a measure of the discrepancy to minimize under those constraints.
Given the hypotheses and the needs, we focus on linear constraints that can be interpreted in term of moments perturbation. This will lead to a quantification of the impact on $P$ of each variable. This will also provide results on the sensitivity of $P$ to the choice of input distributions.
We will later present the perturbations corresponding to a mean shift and a variance shift.
 
Recall that between two pdf $p$ and $q$ we have
\begin{equation}
	KL(p,q)={\displaystyle \int_{-\infty}^{+\infty}p(y)\log\frac{p(y)}{q(y)}dy} \mbox{ if }  \log\frac{p(y)}{q(y)} \in L^1( p(y)dy). 
\end{equation}

Let $i=1,\cdots, d$, the constraints are expressed as follows in function of the
modified density $f_{\mbox{\tiny mod}}$:

\begin{equation}
	\int g_{k}(x_i)f_{\mbox{\tiny mod}}(x_i)dx_i=\delta_{k,i} \; \left(k=1 \cdots K\right)\label{eq:contraintes}.
\end{equation}

Here, for $k=1,\cdots, K$, $g_{k}$ are given functions and $\delta_{k,i}$ are given real.
These quantities will lead to a perturbation of the original density.
The modified density $f_{i\delta}$ considered in our work is:
\begin{equation}
	f_{i\delta}=\mathop{\mbox{argmin}}\limits_{f_{\mbox{\tiny mod}} |(\ref{eq:contraintes}) \mbox{ \tiny{holds}}} KL(f_{\mbox{\tiny mod}},f_i)
\label{eq:problem maxent}
\end{equation}

and the result takes an explicit form (\cite{csiszar1975divergence}) given in the following proposition.

\begin{proposition} 
	\label{2.1}
	Let us define, for $\boldsymbol{\lambda}=\left(\lambda_{1},\cdots,\lambda_{K}\right)^{T}\in\mathbb{R}^K$, 
	\begin{equation}
		\psi_i(\boldsymbol{\lambda})=\log\int f_{i}(x)\exp\left[\sum_{k=1}^{K}\lambda_{k}g_{k}(x)\right]dx\;,  \label{eq:phi}
	\end{equation}
	where the last integral can be finite or infinite (in this last case $\psi_i(\boldsymbol{\lambda})=+\infty$). Further,
	set $\mbox{Dom }\psi_i=\{\boldsymbol{\lambda}\in\mathbb{R}^K | \psi_i(\boldsymbol{\lambda})< +\infty\}$. 
	Assume that there exists at least one pdf $f_m$ satisfying (\ref{eq:contraintes}) 
	and that $\mbox{Dom }\psi_i$ is an open set. Then, there exists a unique $\boldsymbol{\lambda}^*$ such that 
	the solution of the minimisation problem  (\ref{eq:problem maxent}) is 
	\begin{equation}
		\label{eq:regularization function}
		f_{i\delta}(x_i)=f_{i}(x_i)\exp\left[\sum_{k=1}^{K}\lambda^*_{k}g_{k}(x_i)-\psi_i(\boldsymbol{\lambda^*})\right].
	\end{equation}
\end{proposition}

The theoretical technique to compute $\boldsymbol{\lambda}$ is provided in appendix \ref{appLag}. 
Here are presented two kinds of perturbations used further on. 

\paragraph*{Mean shifting}
The first moment is often used to parameterize a distribution. Thus
the first perturbation presented here is a mean shift, that is expressed
with a single constraint:
\begin{equation}
	\int x_if_{\mbox{\tiny mod}}(x_i)dx_i=\delta_i \;. \label{eq:meanshift}
\end{equation}

In term of SA, this perturbation should be used when the user wants
to understand the sensitivity of the inputs to a mean shift - that
is to say ``what if the mean of input $X_{i}$ were $\delta_i$
instead of ${\displaystyle \mathbb{E}}\left[X_{i}\right]$?''.

\begin{proposition}
	\label{prop:expressionmeanshift}
	Considering the constraint (\ref{eq:meanshift}), under the assumptions of Proposition \ref{2.1} the expression of the optimal perturbed density is 
	\begin{equation}
		f_{i\delta_i}(x_i)=\exp(\lambda^* x_i-\psi_i(\lambda^*))f_{i}(x_i)
	\end{equation}
	where $\lambda^*$ is such that equation (\ref{eq:meanshift}) holds. 
\end{proposition}

Note that Equation (\ref{eq:phi}) becomes
\begin{equation}
	\psi_i(\lambda)=\log\int f_{i}(x_i)\exp(\lambda x_i)dx_i=\log\left(M_{X_{i}}(\lambda)\right)
\end{equation}
where $M_{X_{i}}(u)$ is the moment generating function (mgf) of the
$i-$th input. With this notation, $\lambda^*$ is such that
\[
	\int x_i\exp\left(\lambda^* x_i-\log\left(M_{X_{i}}(\lambda^*)\right)\right)f_{i}(x_i)dx_i=\delta_i \;,
\]
which leads to
\[
	\int x_i\exp\left(\lambda^* x_i\right) f_{i}(x_i)dx=\delta_i M_{X_{i}}(\lambda^*) \;.
\]
This can be simplified to:
\begin{equation}
	\frac{M_{X_{i}}'(\lambda^*)}{M_{X_{i}}(\lambda^*)}=\delta_i\;.
\end{equation}
This equation may be easy to solve when the expression of the mgf of the input $X_i$ and of its derivative is known. 

\paragraph*{Variance shifting}
In some cases, the expectation of an input may not be the main source of uncertainty. 
One might be interested in perturbing its second moment.
This case may be treated considering a couple of constraints.
The perturbation presented is a variance shift, therefore the set of constraints is: 
\begin{equation} 
	\begin{cases}
	\int x_if_{\mbox{\tiny mod}}(x_i)dx_i={\displaystyle \mathbb{E}}\left[X_{i}\right]\;,\\
	\int x^{2}_if_{\mbox{\tiny mod}}(x_i)dx_i=V_{\mbox{\tiny  per},i}+{\displaystyle \mathbb{E}}\left[X_{i}\right]^{2}\;.
	\end{cases}\label{eq:cstrvariance}
\end{equation}
The perturbed distribution has the same expectation ${\displaystyle \mathbb{E}}\left[X_{i}\right]$
as the original one and a perturbed variance $V_{\mbox{\tiny per},i}=\mbox{Var} \left[ X_i \right]  \pm \delta_i$.

\begin{proposition} \label{prop:expressionvarshift}
	Under the assumptions of Proposition \ref{2.1},
	for the constraint (\ref{eq:cstrvariance}), the expression of
	the optimal perturbed density is: 
	\[
	f_{i\delta_i}(x_i)=\exp(\lambda^*_{1}x+\lambda^*_{2}x^{2}-\psi_i(\boldsymbol{\lambda^*}))f_{i}(x_i)
	\]
	where $\lambda_{1}^*$ and $\lambda_{2}^*$ are so that equation (\ref{eq:cstrvariance}) holds. 
\end{proposition}

\begin{figure}
\center \includegraphics[width=.75\textwidth]{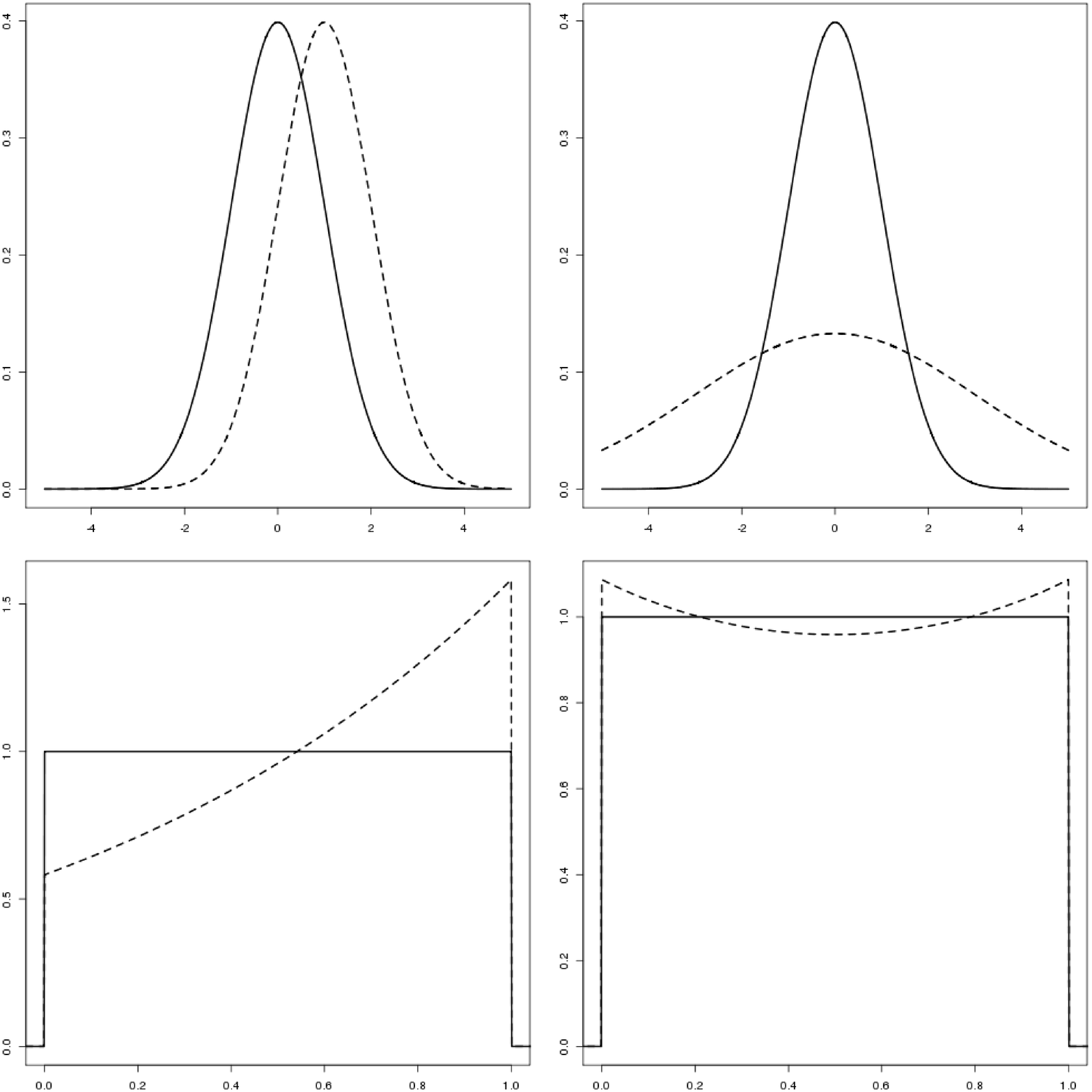} 
\caption{\label{fig:DiversesDensitesModifiees} Mean shifting (left) and variance shifting (right) for Gaussian (upper) and Uniform (lower) distributions. The original distribution is plotted in solid line, the perturbed one is plotted in dashed line. }
\end{figure}

As an example, the two kind of perturbations previously presented are provided for two families of inputs (Gaussian and Uniform) in figure \ref{fig:DiversesDensitesModifiees}.
The perturbations are respectively a mean and variance increasing.
 It is noticeable (and will be proved further on) that the shape is conserved for the Gaussian distribution when shifting the mean or the variance. 
 On the other hand, when increasing its mean, the Uniform distribution is packed down on the right-hand boundary of its support. 
 When increasing its variance, the density is packed down on both boundaries of its support.

\paragraph*{Perturbation of Natural Exponential Family}

In general, when perturbating the input densities, the shape is not conserved. However in the specific case of Natural Exponential Family (NEF), the following proposition can be derived. 

\begin{proposition}\label{prop:NEF}
	Assume that the original random variable $X_i$ belongs to the NEF, i.e. its pdf can be written as:
	\[
	 	f_{i,\theta}(x_i)=b(x_i)\exp\left[x_i\theta-\eta(\theta)\right] \label{eq:notationNEF}	
	\]
	where $\theta$ is a parameter from a parametric space $\Theta$, $b(.)$ is a function that depends only of $x_i$ and 
	\[
		\eta(\theta)=\log\int b(x)\exp\left[x_i\theta\right]dx_i	
	\]
	is the cumulant distribution function.
	Considering the assumptions of Proposition \ref{2.1}, the optimal pdfs proposed respectively in Proposition \ref{prop:expressionmeanshift} and  Proposition \ref{prop:expressionvarshift}
	are also distributed according to a NEF.
\end{proposition}

The proof comes from theorem 3.1 in \cite{csiszar1975divergence}. 
The details of computation are given for a mean shift and a variance shift in Appendix \ref{appedix:proofsofnef}.

% paragraphe sur le cas crue?

%%%%%%%%%%%%%%%%%%%%%%%%%%%%%%%%%%%%%%%%%%%%%%%%%%%%%%%%%%%%%%%%%%%%%%%%%%%%%%%%%%%%%%%%%%%%%%%%%%%%%%%%%%%%%%%
  
\section{Definition, estimation and properties of a sensitivity index\label{sec:Definition}}

Given a unidimensional input variable $X_i$ with pdf $f_i$ and the corresponding perturbed random input $X_{i\delta}\sim f_{i\delta}$.
The perturbed failure probability becomes: 
\begin{equation}
	P_{i\delta}={\displaystyle \int}\1_{\{G(\mathbf{x})<0\}}\dfrac{f_{i\delta}(x_{i})}{f_{i}(x_{i})}f(\mathbf{x})d\mathbf{x}
\end{equation}
where $x_{i}$ is the $i^{\mbox{\begin{scriptsize}th\end{scriptsize}}}$ 
component of the vector $\mathbf{x}$. Independently of the mechanism chosen for the perturbation (see previous section for proposals),
a good sensitivity index $S_{i\delta}$ should have intuitive features that make it appealing to reliability engineers and decision-makers.
We believe that the following definition can fulfil these requirements.  
\begin{definition}
\em{
Define the Density Modification Based Reliability Sensitivity Indices (DMBRSI) the quantity $S_{i\delta}$:

\begin{eqnarray*}
	S_{i\delta} & =& \left[\frac{P_{i\delta}}{P}-1\right]\mathbf{1}_{\{P_{i\delta} \geq P\}}+\left[1-\frac{P}{P_{i\delta}}\right]\mathbf{1}_{\{P_{i\delta}<P\}} 
				\ = \ \frac{P_{i\delta}-P}{P\cdot\mathbf{1}_{\{P_{i\delta} \geq P\}}+P_{i\delta}\cdot\mathbf{1}_{\{P_{i\delta}<P\}}}\;.
\end{eqnarray*}
}
\end{definition}

Firstly, $S_{i\delta}=0$ if $P_{i\delta}=P$, as expected if $X_i$ is a non-influential variable or if $\delta$ expresses a negligible perturbation. 

Secondly, the sign of $S_{i\delta}$ indicates how the perturbation impacts the failure probability qualitatively.
It highlights the situations when $P_{i\delta}>P$ i.e. if the remaining ({\it epistemic})
 uncertainty on the modelling $X_i\sim f_i$ can increase the failure risk.
In this case, the uncertainty on the concerned variable should be more accurately
analysed. Conversely, if if  $P_{i\delta}<P$, $P$ can be interpreted as a conservative assessment of the failure probability, 
with respect to variations of $X_i$. 
In such a case, deeper modelling studies on $X_i$ appear less essential. 

Thirdly, given its sign, the absolute value of $S_{i\delta}$ has simple interpretation and provides a level of the conservatism
 or non-conservatism induced by the perturbation. A value of $\alpha>0$ for the index means that $P_{i\delta}=(1+\alpha)P$.
 If $S_{i\delta}=-\alpha<0$ then $P_{i\delta}=(1/(1+|\alpha|))P$. 

% phrase de transition ? 
The postulated ability of $S_{i\delta}$ to enlighten the sensitivity of $P$ to input perturbations must be tested in concrete cases,
when  an estimator $\hat{P}_N$ of $P$ can be computed using an already available design of $N$ numerical experiments.
In this paper, $N$ is assumed to be large enough such that statistical estimation stands within the framework of asymptotic theory.
Besides, we assume for simplicity a standard Monte Carlo design of experiments, according to which 
$\hat{P}_{N}  =  \sum_{n=1}^{N}\1_{\{G(\mathbf{x}^{n})<0\}}/N$
where the $\mathbf{x}^{1},\cdots,\mathbf{x}^{N}$ are independent realisations of $X$.
 The strong Law of Large Numbers (LLN) and the Central Limit Theorem (CLT) ensure that for almost all realisations $\hat{P}_{N}\xrightarrow[N\rightarrow\infty]{} P$ and  
\begin{eqnarray*}
	\sqrt{N/[P(1-P)]}(\hat{P}_{N}-P) & \xrightarrow[N\rightarrow\infty]{\mathcal{L}} & \mathcal{N}(0,1). 
\end{eqnarray*}
The Monte Carlo framework allows $P_{i\delta}$ to be consistently estimated without new calls to $G$, through a ``reverse'' importance sampling mechanism:
\begin{eqnarray*}
	\hat{P}_{i\delta N} & = & {\frac{1}{N}\sum_{n=1}^{N}{\1}_{\{G(\mathbf{x}^{n})<0\}}\frac{f_{i\delta}(x_{i}^{n})}{f_{i}(x_{i}^{n})}}.
\end{eqnarray*}
This property holds in the more general case when $P$ is originally estimated by importance sampling rather than simple Monte Carlo, which is more appealing when $G$ is time-consuming \cite{BeckmanetMcKey87,hesterberg1996estimates}. This generalization is discussed further in the text (Section \ref{sec:Discussion}). 
 The following lemma ensures the asymtotic behaviour of such an estimator.

\begin{lemma}\label{lemma1}
	Assume the usual conditions
	\begin{enumerate}
	\item[(i)] $\mbox{Supp}(f_{i\delta})\subseteq\mbox{Supp}(f_{i})$,
	\item[(ii)] ${\displaystyle \int_{{\Supp}(f_{i})} \frac{f^2_{i\delta}(x)}{f_{i}(x)} \ dx < \infty}$,
	\end{enumerate}
	then $\hat{P}_{i\delta N} \xrightarrow[N\rightarrow\infty]{} P_{i\delta}$ and $\sqrt{N}{\sigma}^{-1}_{i\delta N}\left(\hat{P}_{i\delta N}- P_{i\delta}\right)  \xrightarrow[N\rightarrow\infty]{\mathcal{L}}   {\cal{N}}(0,1).$
	The exact expression of ${\sigma}^{-1}_{i\delta N}$ is given in Appendix \ref{proofs}, equation (\ref{eq:expression_sigma}). It can be consistently estimated by 
	\begin{eqnarray*}
	\hat{\sigma}^{2}_{i\delta N} & = & \frac{1}{N}\sum\limits_{n=1}^N \1_{\{G(\mathbf{x}^{n})<0\}} \left(\frac{f_{i\delta}(x_{i}^{n})}{f_{i}(x_{i}^{n})}\right)^2  - \hat{P}^2_{i\delta N}. \\
	\end{eqnarray*}
	The proof of this Lemma is given in Appendix \ref{proofoflemma1}.
\end{lemma}

\noindent The asymptotic properties of any estimator of $S_{i\delta}$ will depend on the correlation between $\hat{P}_N$ and $\hat{P}_{i\delta N}$. The next proposition summarizes the features of the joint asymptotic distribution of both estimators.

\begin{proposition}\label{CTL.correlated}
	Under assumptions (i) and (ii) of Lemma \ref{lemma1},
	\begin{eqnarray*}
	\sqrt{N}\left[\left(\begin{array}{c}
		\hat{P}_{N}\\
		\hat{P}_{i\delta N}
	\end{array}\right)-\left(\begin{array}{c}
		P\\
		P_{i\delta}
	\end{array}\right)\right]  \xrightarrow[N\rightarrow\infty]{\mathcal{L}}  \mathcal{N}_{2}\left({0},\Sigma_{i\delta}\right)
	\end{eqnarray*}
	where $\Sigma_{i\delta}$ is given in Appendix \ref{proofs}, Equation (\ref{eq:expression_Sigma}) and can be consistently estimated by
	\begin{eqnarray*}
		\hat{\Sigma}_{i\delta} & = & \left(\begin{array}{cc}
		\hat{P}_{N}(1-\hat{P}_{N}) & \hat{P}_{i\delta N}(1-\hat{P}_{N})\\
		\hat{P}_{i\delta N}(1-\hat{P}_{N}) & \hat{\sigma}^{2}_{i\delta N}
		\end{array}\right). \\
	\end{eqnarray*}
	The proof of this Proposition is given in Appendix \ref{A2}.
\end{proposition} 

\noindent Given $(\hat{P}_{N},\hat{P}_{i\delta N})$, the plugging estimator for $S_{i\delta}$ is
\begin{eqnarray}
	\hat{S}_{\text{i\ensuremath{\delta}}N} & = & \left[\frac{\hat{P}_{i\delta N}}{\hat{P}_{N}}-1\right]\mathbf{1}_{\left\{ \hat{P}_{i\delta N} \geq \hat{P}_{N}\right\} 	}+\left[1-\frac{\hat{P}_{N}}{\hat{P}_{i\delta N}}\right]\mathbf{1}_{\left\{ \hat{P}_{i\delta N}<\hat{P}_{N}\right\} }. \label{S.estimator}
\end{eqnarray}
In corollary of Proposition \ref{CTL.correlated}, applying the continuous-mapping theorem to the  function 
$s(x,y)  =  \left[\frac{y}{x}-1\right]\mathbf{1}_{y \geq x}+\left[1-\frac{x}{y}\right]\mathbf{1}_{y<x},$
$\hat{S}_{\text{i\ensuremath{\delta}}N}$ converges almost surely to $S_{i\delta}$. 

The following CLT results from Theorem 3.1 in \cite{van2000asymptotic}. 
\begin{proposition} 
	Assume that assumptions (i) and (ii) of Lemma \ref{lemma1} hold and further that $P\neq P_{i\delta}$, we have
	\begin{equation}
		\sqrt{N}\left[\hat{S}_{\text{i\ensuremath{\delta}}N}-S_{\text{i\ensuremath{\delta}}}\right]\xrightarrow[N\rightarrow\infty]{\mathcal{L}}\mathcal{N}\left(0,d^{T}\Sigma d\right)
	\end{equation}
	with $\displaystyle d=\left(\frac{\partial s}{\partial x}(P,P_{i\delta}),\frac{\partial s}{\partial y}(P,P_{i\delta})\right)^{T}$ for $x\neq y$, and
	\begin{eqnarray*}
		\frac{\partial s}{\partial x}(x,y) & = & -y\1_{\{y \geq x\}}/{x^{2}} - \frac{1}{y}} \1_{\{y<x\},\\
		\frac{\partial s}{\partial y}(x,y) & = & \frac{1}{x}\1_{\{y \geq x\}} + x\1_{\{y<x\}}/{y^{2}}. \\
	\end{eqnarray*}
	This holds when $P= P_{i\delta}$. Indeed, one has for $x^* \neq 0$ :
	$$
	\lim_{\pileg} \nabla s(x,y)  =  \lim_{\piled} \nabla s(x,y) = \left(-\frac{1}{x^*} , \frac{1}{x^*} \right)^T.
	$$
\end{proposition}

% § sur le cas crue? peut ^etre pas la peine dans cette section.

%%%%%%%%%%%%%%%%%%%%%%%%%%%%%%%%%%%%%%%%%%%%%%%%%%%%%%%%%%%%%%%%%%%%%%%%%%%%%%%%%%%%%%%%%%%%%%%%%%%%%%%%%%%%%%%

\section{Numerical experiments\label{sec:Numerical-experiments}}

In this section, the methodology presented throughout the article is tested on two academic cases and a more realistic industrial model.
 The new indices are compared to the results provided by two reference methods, FORM indices (or importance factors)
and Sobol' indices. First, we  briefly present these two methods.

\subsection{Presentation of the two reference SA methods}

\subsubsection{FORM importance factors}
The First Order Reliability Method (FORM)\cite{lemaire2009structural} is an estimation technique for a failure probability based upon approximating the failure domain and solving an optimization problem.
In practice, it is considered a standard solution in structural reliability, mainly because of its low cost.
This method proceeds in four steps:
\begin{itemize}
	\item transformation of the input variables space into the standard space (a space where all the variables are standard independent Gaussian);
	\item search of the closest failure point to the origin of the standard space (also referred as the design point) via an optimisation algorithm;
	\item approximation of the failure surface by an hyperplane;
	\item failure probability estimation based on the geometry of the failure domain.
\end{itemize}
Importance factors are byproducts of this method. These sensitivity measures aims at quantifying the importance of a variable on the failure probability.
From the design point $\mathbf{u}^{*}$ one writes:
\begin{equation*}
	\mathbf{u}^{*}=\beta_{HL}\mathbf{\mbox{$\alpha$}}^{*}
\end{equation*} 
where $\beta_{HL}$ is the distance between the origin of the standard space and $\mathbf{u}^{*}$ and $\alpha^{*}$
is the normalised vector of direction. Then for each variable  $U_{i}$, one can obtain the importance factor $\alpha_{i}^{*^{2}}$. 

\subsubsection{Sobol' sensitivity indices}
In the SA framework, let us have $\mathbf{X}=(X_{1,}...,X_{d})$, a random vector where the variables are mutually independent
 and $Y=G(\mathbf{X})$, the output of a deterministic code $G()$. Thus a functional decomposition of the variance is feasible, often 
 referred as functional ANOVA: $\mbox{Var}[Y]=\sum_{i=1}^{d}V_{i}(Y)+\sum_{i<j}^{d}V_{ij}(Y)+\cdots+V_{12...d}(Y)$
 where $V_{i}(Y)=\mbox{Var}[\mathbb{E}(Y|X_{i})]$, $V_{ij}(Y)=\mbox{Var}[\mathbb{E}(Y|X_{i},X_{j})]-V_{i}(Y)-V_{j}(Y)$
  and so on for higher order interactions.
The so-called “Sobol' indices” or “sensitivity indices” \cite{sobol1993sensitivity} are obtained as follows:
\begin{equation*}
	S_{i}=\frac{V_{i}(Y)}{\mbox{Var}[Y]},\quad S_{ij}=\frac{V_{ij}(Y)}{\mbox{Var}[Y]},\quad\cdots
\end{equation*}
These indices express the share of variance of $Y$ that is due to a given input or a set of inputs.
The number of indices growths in an exponential way with the number $d$ of dimension: there are $2^{d}-1$ indices.
For computational time and interpretation reasons, the practitioner rarely estimate indices of order higher than $2$.
Therefore Homma \& Saltelli \cite{homma1996importance} introduced the so-called “total Sobol index” or “total effect” that writes as follow:
\begin{equation*}
	S_{T_{i}}=S_{i}+\sum_{j\neq i}S_{ij}+\sum_{j\neq i,k\neq i,j<k}S_{ijk}+...
\end{equation*} 
This total sensitivity index contains all the effects due to $X_i$ and the interactions between $X_i$ and other inputs. 

\subsubsection{Testing methodology}
Importance Factors and Sobol' indices (SI)  are computed using the methodologies given in \cite{lemaire2009structural}
and \cite{sobol1993sensitivity}, respectively. The Sobol' indices are computed using two initial samples of size $10^5$,
 resulting into $N=10^5\times(d+2)$ function calls (\cite{saltelli2010variance}).
To assess the reproductibility of the estimation of the SI, $50$ replications are made.
The averages of the obtained values and the coefficients of variation (C.o.v.) of the indices are provided.
One should notice that the SI are applied on the indicator of the failure function $\1_{\{G(\mathbf{x})<0\}}$.
Following the definitions of the Importance Factors and SI, those indices lie in $\left[0,1\right]$.

\subsection{Hyperplane failure surface}

For the first example, $\mathbf{X}$ is set to be a
$4-$dimensional vector, with $d=4$ independent marginal distributions
normally distributed with parameters $0$ and $1$. Therefore $f_{X_{i}}\sim\mathcal{N}(0,1)$ for $i=1,..,4$. The failure function is defined as:
$$G(\mathbf{X})=k-\sum_{i=1}^{4}a_{i}X_{i}$$ where $k$ and $\mathbf{a}=(a_{1},a_{2},a_{3},a_{4})$
are the parameters of the model. For this numerical example, parameters
are set with values $k=16$ and $\mathbf{a}=(1,-6,4,0)$. An explicit
expression for $P$ can be given:
$$
P=\phi\left(-k/{\displaystyle \sqrt{\sum_{i=1}^{4}a_{i}^{2}}}\right)=\phi\left(\frac{-16}{{\displaystyle \sqrt{53}}}\right)\simeq0.014
$$
where $\phi(.)$ is the standard normal cumulative distribution function.
  
It is expected that the influence of $X_i$ on $P$ uniquely depends on $|a_i|$.
 The larger the absolute value of the coefficient, the larger the expected influence. 
 The aim for choosing one non-influential variable $X_4$ (because $a_4=0$) is to assess if the SA methods can identify this variable as non-influential on the failure probability.

\subsubsection{FORM }
%The algorithm FORM is launched starting from the origin point of the Gaussian space. It converges to the true and unique design
%point $\left(0.30,-1.81,1.2,0\right)$ in $39$ function calls, giving
%an approximate probability of $\hat{P}_{FORM}=0.01398$. The importance
%factors are displayed in table \ref{tab:IFhyper}. The very good performance
%of FORM to estimate the failure probability can be noted, which is
%expected given that the failure surface is an hyperplane.

In this ideal hyperplane failure surface case, FORM provides an approximated value $\hat{P}_{FORM}=0.01398$, which is as expected (\cite{lemaire2009structural}) close to the exact value. 
$39$ model calls have been required. 
The importance factors, given in Table \ref{tab:IFhyper}, provide an accurate variable ranking for the failure function.

{\footnotesize }
\begin{table}[h]
\begin{centering}
{\footnotesize }%
\begin{tabular}{|c|c|c|c|c|}
\hline 
{\footnotesize Variable} & {\footnotesize $X_{1}$} & {\footnotesize $X_{2}$} & {\footnotesize $X_{3}$} & {\footnotesize $X_{4}$}\tabularnewline
\hline 
\hline 
{\footnotesize Importance factor} & {\footnotesize $0.018$} & {\footnotesize $0.679$} & {\footnotesize $0.302$} & {\footnotesize $0$}\tabularnewline
\hline 
\end{tabular}
\par\end{centering}{\footnotesize \par}

{\footnotesize \caption{\label{tab:IFhyper}{\footnotesize Importance factors for hyperplane
function}}
}
\end{table}
{\footnotesize \par}

\subsubsection{Sobol' indices}
The first-order and total indices are displayed in Table \ref{tab:SobolHyper}.  The interpretation of the results is that  $X_{2}$ and $X_{3}$ concentrate most of the variance of the indicator function. At first order, $25\%$ of its variance is explained by $X_{2}$ without any interaction. It should be noted that the total index for $X_{4}$ is null, assessing that this variable does not impact the failure probability. The C.o.v. of the total indices estimators are small, meaning that this method is reproducible and that $6\times10^{5}$ points are enough to estimate in an efficient way the indices $S_{Ti}$. On the other hand, some C.o.v. values for small first order indices are quite high. This effect is well-known in SI estimation problems and can be corrected by improved formulas \cite{saltelli2010variance}.

% discuter le fait que les indices totaux sont moins discriminant que les facteurs form

%The conclusion of this result is that the method correctly estimates high indices but estimates poorly the indices close to $0$. On the other hand, the relevant information is that the index is close to $0$. Thus this situation may not be a problem. 

{\footnotesize }
\begin{table}[h]
\label{tab:SobolHyper}
\begin{centering}
{\footnotesize }%
\begin{tabular}{|c|c|c|c|c|c|c|c|c|}
\hline 
{\footnotesize Sobol' Index} & {\footnotesize $S_{1}$} & {\footnotesize $S_{2}$} & {\footnotesize $S_{3}$} & {\footnotesize $S_{4}$} & {\footnotesize $S_{T1}$} & {\footnotesize $S_{T2}$} & {\footnotesize $S_{T3}$} & {\footnotesize $S_{T4}$}\tabularnewline
\hline 
\hline 
{\footnotesize Mean} & {\footnotesize $0.0017$} & {\footnotesize $0.2575$} & {\footnotesize $0.0544$} & {\footnotesize $9.45.10^{-5}$} & {\footnotesize $0.1984$} & {\footnotesize $0.9397$} & {\footnotesize $0.7256$} & {\footnotesize $0$}\tabularnewline
\hline 
{\footnotesize C.o.v.} & {\footnotesize $1.5854$} & {\footnotesize $0.04826$} & {\footnotesize $0.1336$} & {\footnotesize $27.4$} & {\footnotesize $0.012$} & {\footnotesize $0.0069$} & {\footnotesize $0.013$} & {\footnotesize $0$}\tabularnewline
\hline 
\end{tabular}
\par\end{centering}{\footnotesize \par}

{\footnotesize \caption{{\footnotesize Sobol' indices for hyperplane function}}
}
\end{table}
{\footnotesize \par}

\subsubsection{Density modification based reliability indices}
The method presented throughout this article is applied on the hyperplane function. As explained in section \ref{sec:Methodologies-of-input}, several ways to perturb the input distributions exist. For this case, a mean shifting is first applied, then a variance shifting with fixed mean. A simple calculus gives that the perturbed pdf are Gaussian, respectively with the constraint mean and variance 1 for the mean shifting perturbation (see Table \ref{tab:TableRecap}), and with mean 0 and the constraint variance for the variance shifting perturbation. 

\begin{figure}
\includegraphics[width=1\textwidth]{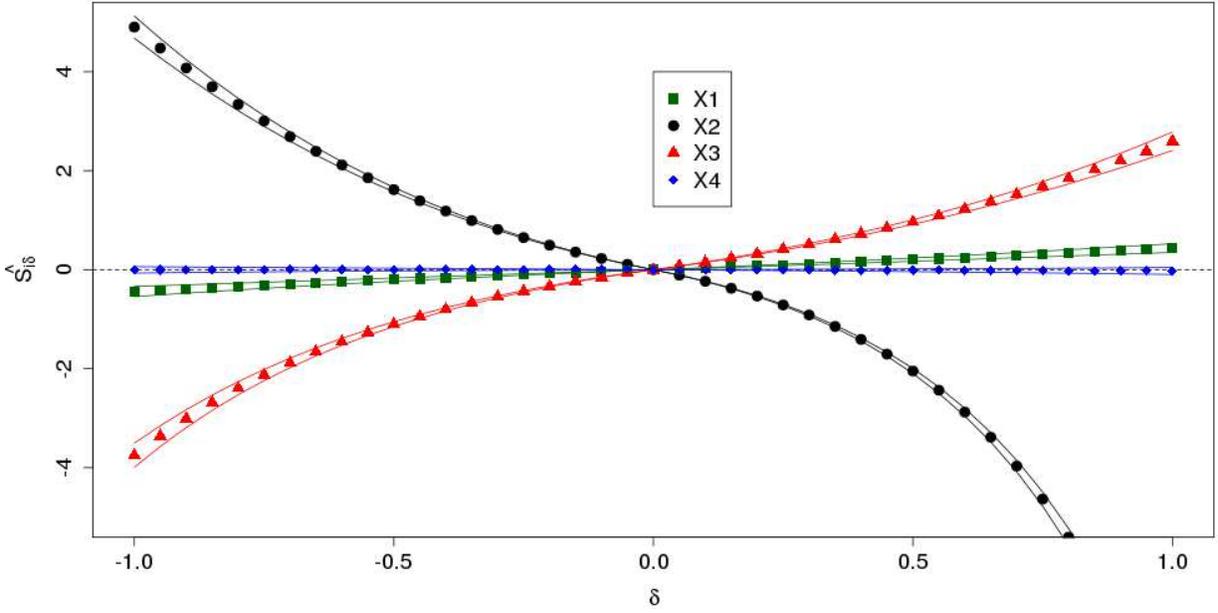} 
\caption{\label{fig:IBSMDHyperMean} Estimated indices $\widehat{S_{i\delta}}$ for hyperplane function with a mean shifting}
\end{figure}

\begin{figure}
\includegraphics[width=1\textwidth]{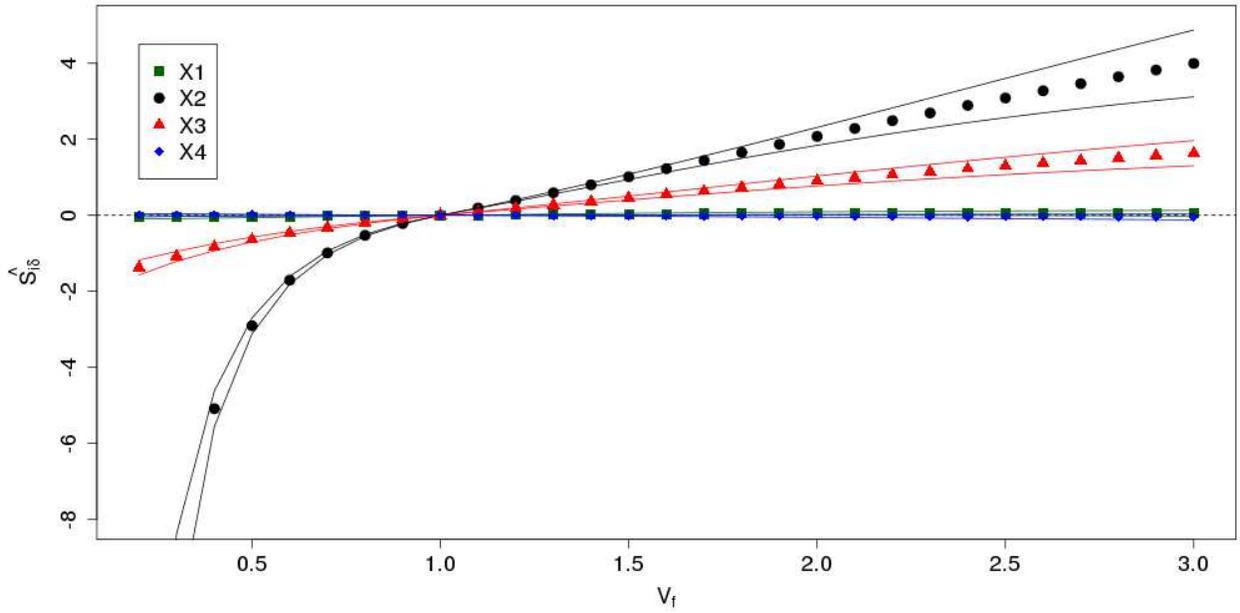}
\caption{\label{fig:IBSMDHyperVariance}
Estimated indices $\widehat{S_{i,V_{\mbox{f}}}}$ for hyperplane function with a variance shifting}
\end{figure}

The MC estimation gives $\hat{P}=0.01446$ with $10^5$ function calls. For the mean shifting (see (\ref{eq:meanshift})), the domain variation for $\delta$ ranges from $-1$ to $1$ with $40$ points, reminding that $\delta=0$ cannot be considered as a perturbation. For the variance shifting (see (\ref{eq:cstrvariance})), the variation domain for $V_{\mbox{\tiny per}}$ ranges from $1/20$ to $3$ with $28$ points, where  $V_{\mbox{\tiny per}}=1$ is not a perturbation. The estimated indices are plotted respectively in Figure \ref{fig:IBSMDHyperMean} for mean shifting and in Figure \ref{fig:IBSMDHyperVariance} for variance shifting. The 95\% confidence intervals are plotted around the indices, using the presented asymptotic formulas in Section \ref{sec:Definition}.

\paragraph*{Mean perturbation indices}
The indices $\widehat{S_{i\delta}}$ behave in a monotonic way given the importance of the perturbation. The slope at the origin is directly related to
the value of $a_{i}$. For influential variables ($X_{2}$ and $X_{3}$), the increasing or the decreasing is faster than linear, whereas the curve seems linear for the slightly influential variable ($X_{1}$). 
Modifying the mean with a positive amplitude slightly rises the failure probability for  $X_{1}$, highly decreases it for  $X_{2}$ and increases it for $X_{3}$ (Figure \ref{fig:IBSMDHyperMean}). The effects are reversed with similar amplitude for negative $\delta$. It can be seen that $X_{4}$ has no impact on the failure probability for any perturbation. Those results are consistent with the expression of the failure function.
One can see that the confidence intervals (CI) associated to $X_2$ and $X_3$ are fairly well separated, except for the small absolute value of $\delta$. On the other hand, the CI associated to $X_1$ and $X_4$ are not separated until absolute value of $\delta$ higher than $0.6$. It can be concluded from the observation of the CI is that the impact of variable $X_1$ and $X_4$ cannot be differentiated, unless a broad change of the mean occurs. 
\paragraph*{Variance perturbation indices}
Increasing the variance of input $X_{2}$ and $X_{3}$ increases the failure probability, whereas it decreases when decreasing the variance (Figure \ref{fig:IBSMDHyperVariance}). Modifying the variance of $X_{1}$ and $X_{4}$ have no effect on the failure probability. The increasing of the indices is linear for $X_{2}$ and $X_{3}$, and the decreasing of the indices is faster than linear, especially for $X_{2}$.  
Considering the CI, one can see that they are well separated for variable $X_2$ and $X_3$, assessing the relative importance of these variables. On the other hand, the CI associated to $X_1$ and $X_4$ are not separated and contain $0$.

\subsubsection{Conclusion}

The DMBRSI has brought the following conclusions: when shifting the mean (that is to say the central tendency in this case), the most influential variable is $X_2$, followed by $X_3$. $X_1$ is slightly influential while $X_4$ is not influential at all. When shifting the variance, variable $X_2$ is more influential than variable $X_3$. Variables $X_1$ and $X_4$ have no impact when shifting the variance.
We argue that all these information are much richer than the ones provided by importance factors and by Sobol' indices. Indeed, the information are provided about regions of the input space leading to failure event.
This is, in our opinion, more of interest to the practitioner than a "simple" variable ranking.

\subsection{Thresholded Ishigami function}
 
\begin{figure}
\includegraphics[width=1\textwidth]{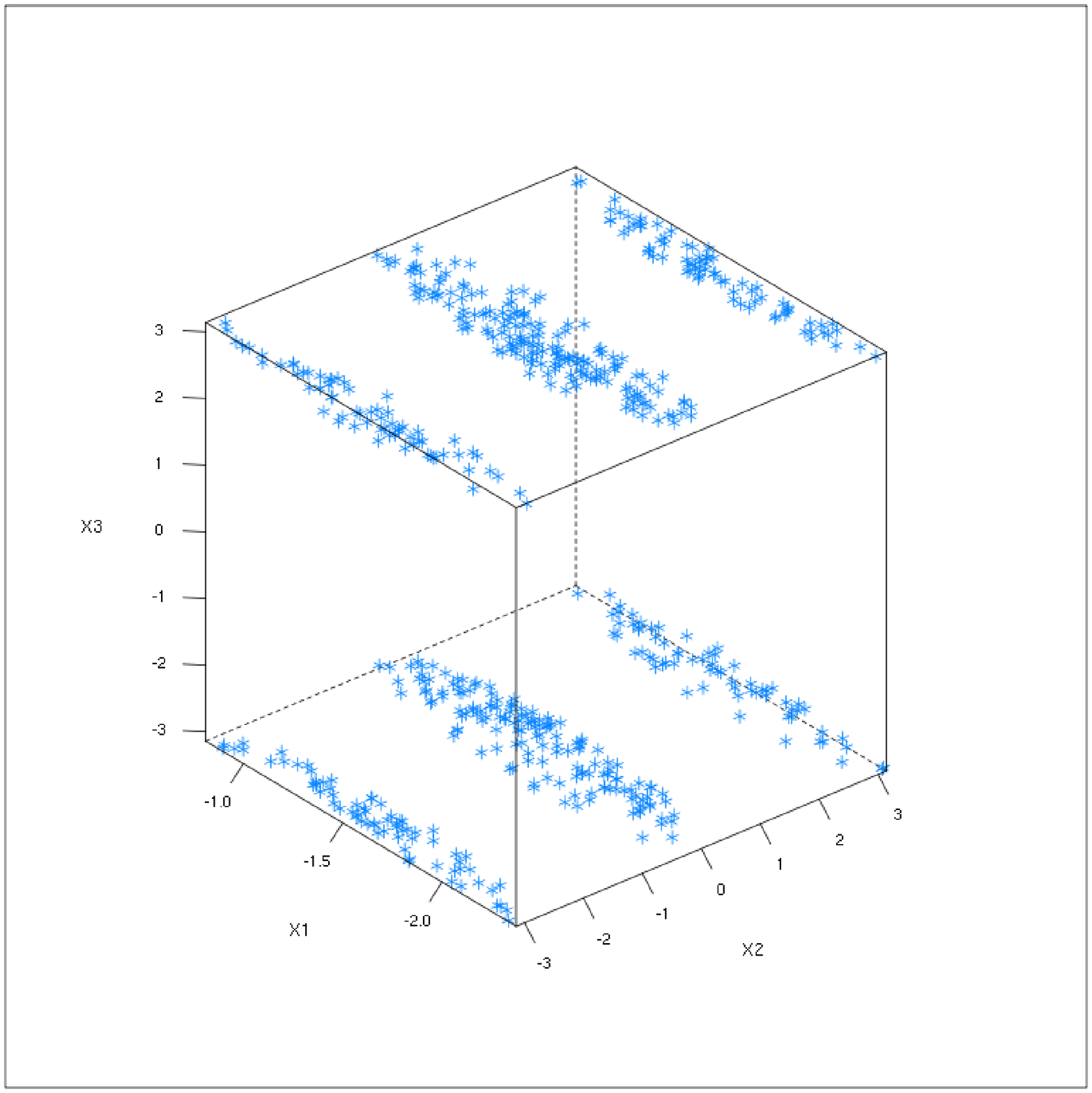}
\caption{\label{fig:ishidots} Ishigami failure points from a MC sample}
\end{figure}

The Ishigami function is a common test case in SA since it presents high degree of non-linearity with interactions between the variables. A modified version of the Ishigami function will be considered in this paper. One has
$$
G(\mathbf{X})=\sin\left(X_{1}\right)+7\sin\left(X_{2}\right)^{2}+0.1X_{3}^{4}\sin\left(X_{1}\right)+7
$$
where $\mathbf{X}$ is a $3-$dimensional vector of independent marginals uniformly distributed on $\left[-\pi,\pi\right].$ In Figure \ref{fig:ishidots}, the failure points (where $G(\mathbf{x})<0$) are plotted in a 3-d scatterplot.

There are $614$ failure points on a MC sample of $10^{5}$ points. Therefore the failure probability here is roughly $\hat{P}=6.14\times10^{-3}$. The complex repartition of the failure points can be noticed. Those points lay in a zone defined by the negative values of $X_{1}$, the extremal and mean values of $X_{2}$ (around $-\pi$, $0$ and $\pi$), and the extremal values of $X_{3}$ (around $-\pi$ and $\pi$).

\subsubsection{FORM} 
The algorithm FORM converges to an incoherent design point $\left(6.03,0.1,0\right)$ in $50$ function calls, giving an approximate probability of $\hat{P}_{FORM}=0.54$. The importance factors are displayed in Table \ref{tab:IFIshi}. The  bad performance of FORM is expected given that the failure domain consists in six separate domains and that the function is highly oscillant, leading to optimization difficulties. The design point is aberrant, therefore the FORM results of SA are incorrect.

{\footnotesize }
\begin{table}[h]
\begin{centering}
{\footnotesize }%
\begin{tabular}{|c|c|c|c|}
\hline 
{\footnotesize Variable} & {\footnotesize $X_{1}$} & {\footnotesize $X_{2}$} & {\footnotesize $X_{3}$}\tabularnewline
\hline 
\hline 
{\footnotesize Importance factor} & {\footnotesize $1e^{-17}$} & {\footnotesize $1$} & {\footnotesize $0$}\tabularnewline
\hline 
\end{tabular}
\par\end{centering}{\footnotesize \par}
{\footnotesize \caption{\label{tab:IFIshi}{\footnotesize Importance factors for Ishigami
function}}
}
\end{table}
{\footnotesize \par}

\subsubsection{Sobol' indices} 
The first-order and total indices, computed with $5\times10^5$ function calls are displayed in Table \ref{tab:SobolIshi}. The small values of first order indices show that no variable has impact on the variance of the indicator of failure on its own. The three total indices have relatively high and similar values. This states that all the variables highly interact with each other to cause system failure. The SI method is thus non-discriminant in this case. The low C.o.v. show that the method is reproducible.

{\footnotesize }
\begin{table}[h]
\begin{centering}
{\footnotesize }%
\begin{tabular}{|c|c|c|c|c|c|c|}
\hline 
{\footnotesize Sobol' Index} & {\footnotesize $S_{1}$} & {\footnotesize $S_{2}$} & {\footnotesize $S_{3}$} & {\footnotesize $S_{T1}$} & {\footnotesize $S_{T2}$} & {\footnotesize $S_{T3}$}\tabularnewline
\hline 
\hline 
{\footnotesize Mean} & {\footnotesize $0.0234$} & {\footnotesize $0.0099$} & {\footnotesize $0.0667$} & {\footnotesize $0.8158$} & {\footnotesize $0.6758$} & {\footnotesize $0.9299$}\tabularnewline
\hline 
{\footnotesize C.o.v.} & {\footnotesize $0.0072$} & {\footnotesize $0.0051$} & {\footnotesize $0.0095$} & {\footnotesize $0.0156$} & {\footnotesize $0.0216$} & {\footnotesize $0.0094$}\tabularnewline
\hline 
\end{tabular}
\par\end{centering}{\footnotesize \par}

{\footnotesize \caption{\label{tab:SobolIshi}{\footnotesize Sobol' indices for Ishigami function}}
}
\end{table}
{\footnotesize \par}

\subsubsection{Density modification based reliability indices}
\begin{figure}
\includegraphics[width=1\textwidth]{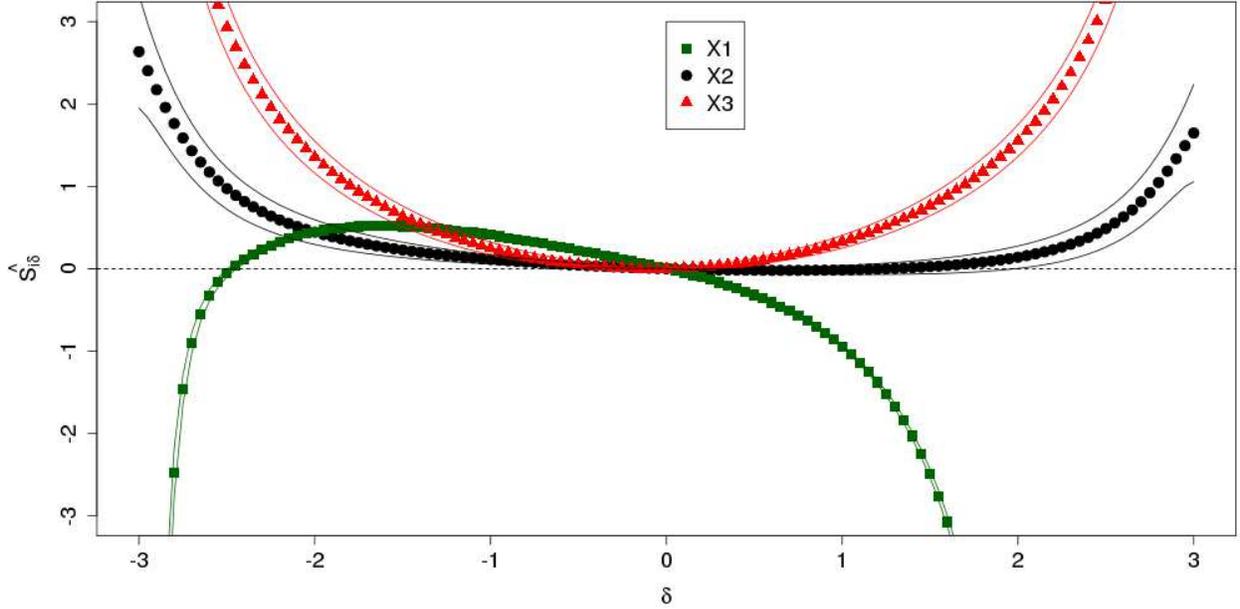}
\caption{\label{fig:IBSMDIshiMean}
Estimated indices $\widehat{S_{i\delta}}$ for the thresholded Ishigami function with a mean shifting}
\end{figure}

\begin{figure}
\includegraphics[width=1\textwidth]{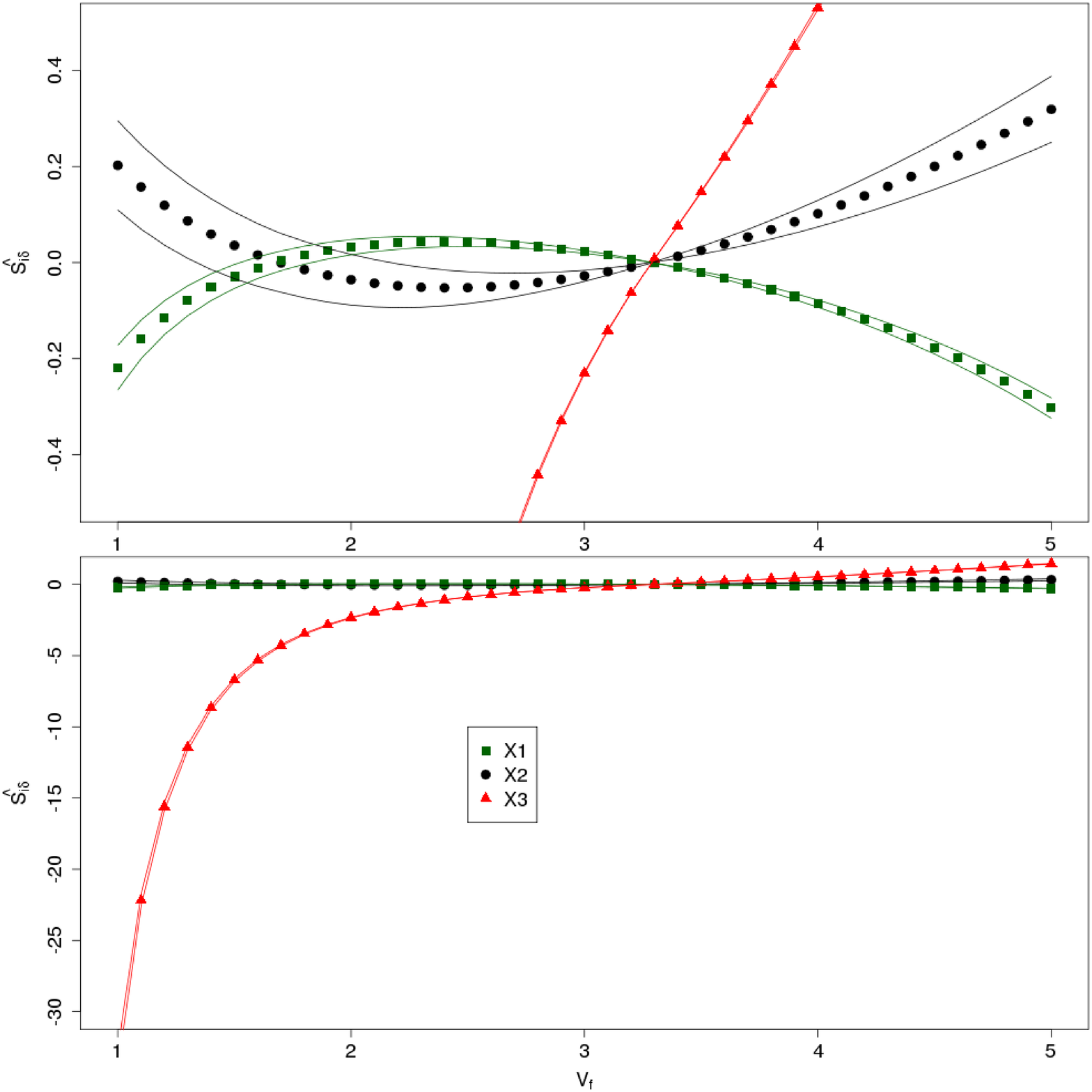}
\caption{\label{fig:IBSMDIshiVariance}
Estimated indices $\widehat{S_{i,V_{\mbox{\tiny per}}}}$ for the  thresholded Ishigami function with a variance shifting. The upper figure is a zoom where the $\widehat{S_{i,V_{\mbox{\tiny per}}}}$ axis lies into $\left[-0.5,0.5\right]$. The lower figure shows almost the whole range variation for $\widehat{S_{i,V_{\mbox{\tiny per}}}}$. The curves cross for the value of $V_{\mbox{\tiny per}}$ that corresponds to the original variance, namely $\pi^2/2$}
\end{figure}

 The method presented throughout this article is applied on the thresholded Ishigami function.
  As for the hyperplane test case, a mean shifting and a variance shifting are applied.
  The modified distribution when a mean shift is applied on a uniform distribution is given in Table \ref{tab:TableRecap}.
   The modified pdf when shifting the variance and keeping the same expectation is proportional to a truncated Gaussian when decreasing the variance. When increasing the variance, the perturbed distribution is a symmetrical distribution with 2 modes close to the endpoints of the support (see figure \ref{fig:DiversesDensitesModifiees}).
    As previously, the same MC sample of size $10^{5}$ (also used to produce Figure \ref{fig:ishidots}) is used to estimate the indices with both perturbations.
     For the mean shifting (see (\ref{eq:meanshift})), the variation domain for $\delta$ ranges from  $-3$ to $3$ with $60$ points - numerical consideration forbidding to choose a shifted mean closer to the endpoints. 
     For variance shifting, the variation domain for $V_{\mbox{\tiny f}}$ ranges from $1$ to $5$ with $40$ points. Let us recall that the original variance is $\mbox{Var}[{X_i}]=\pi^{2}/3\simeq3.29.$
      The estimated indices are plotted respectively in Figure \ref{fig:IBSMDIshiMean} for mean shifting and in Figure \ref{fig:IBSMDIshiVariance} for variance shifting.

\paragraph*{Mean perturbation indices}
 A perturbation of the mean for $X_{2}$ and $X_{3}$ will increase the failure probability, though the impact for the same mean perturbation is stronger for $X_{3}$ ($\widehat{S_{3,-3}}$ and $\widehat{S_{3,3}}$ approximately equal respectively $9.5$ and $10$, Figure \ref{fig:IBSMDIshiMean}). On the other hand, the indices concerning $X_{1}$ show that a mean shift between $-1$ and $-2$ increases the failure probability, whereas an increasing of the mean or a large decreasing strongly diminishes the failure probability ($\widehat{S_{1,3}}$ approximatively equals $-7.10^{11}$). Therefore, Figure \ref{fig:IBSMDIshiMean} leads to two conclusions. First, the failure probability can be strongly reduced when increasing the mean of the first variable  $X_{1}$ (this is also provided by Figure \ref{fig:ishidots} wherein all failure points have a negative value of $X_{1}$). Second, any change in the mean for $X_{2}$ or $X_{3}$ will lead to an increase of the failure probability. 
 The CI are well separated, except in the $-1$ to $1$ zone. One can notice that the CI associated to $X_2$ contains 0 between values of $\delta$ from $-1.5$ to $1.5$, thus the associated indices might be null in these case. This has to be taken into account when assessing the relative importance of $X_2$.  

\paragraph*{Variance perturbation indices}
Figure \ref{fig:IBSMDIshiVariance} (upper part) shows that a change in the variance has little effect on  $X_{2}$ and $X_{1}$, though the change is of opposite  effect on the failure probability. However, considering that the indices $\widehat{S_{2,V_{\mbox{\tiny per},i}}}$ and $\widehat{S_{1,V_{\mbox{\tiny per},i}}}$ lie between $-0.4$ and $0.4$, one can conclude that the variance of theses variables are not of great influence on the failure probability. On the other hand, Figure \ref{fig:IBSMDIshiVariance} (lower part) shows that any reduction of $\mbox{Var} \left[ X_{3} \right]$ strongly decreases the failure probability, and that an increase of the variance slightly increases the failure probability. This is relevant with the expression of the failure surface, as $X_{3}$ is fourth powered and multiplied by the sinus of $X_{1}$. A variance decreasing as formulated gives a distribution concentrated around 0. Decreasing  $\mbox{Var} \left[ X_{3} \right]$ shrinks the concerned term in $G(\mathbf{X})$. Therefore it reduces the failure probability. 
The CI associated to $X_3$ are broadly separated from the others.

\subsection{\label{sbuseq:flood}Industrial case : flood case}
The goal of this test case is to assess the risk of a flood over a dyke for the safety of industrial installations \cite{bernardara2010uncertainty}. This comes down to model the level of a flood. As a function of hydraulical parameters, many of them being randomized to account for uncertainty. From a simplification of the Saint-Venant equation, a flood risk model is obtained. The quantity of interest is the difference between the level of the dyke and the height of water. If this quantity is negative, the installation is flooded.
Hydraulical parameters are the following: $Q$ the flow rate, $L$ the watercourse section length studied, $B$ the watercourse width, $K_{s}$ the watercourse bed friction coefficient (also called Strickler coefficient), $Z_{m}$ and $Z_{v}$ respectively the upstream and downstream bottom watercourse level above sea level and $H_{d}$ the dyke height measured from the bottom of the watercourse bed. The water level model is expressed as:  
$$
H=\left(\frac{Q}{K_{s}B\sqrt{\frac{Z_{m}-Z_{v}}{L}}}\right)^{\frac{3}{5}}.
$$
Therefore the following quantity is considered: 
$$G=H_{d}-(Z_{v}+H).$$
Among the model inputs, the choice is made that the
following variables are known precisely: $L=5000$ (m), $B=300$ (m), $H_{d}=58$ (m), and the following are considered to be random. $Q$ ($\text{m}{}^{3}.\text{s}^{-1}$) follows a positively truncated Gumbel distribution of parameters $a=1013$ and $b=558$ with a minimum value of $0$. $K_{s}$ ($\text{m}{}^{1/3} \text{s}^{-1}$) follows a truncated Gaussian distribution of parameters $\mu=30$ and $\sigma=7.5$, with a minimum value of $1$. $Z_{v}$ (m) follows a triangular distribution with minimum $49$, mode $50$ and maximum $51$. $Z_{m}$ (m) follows a triangular distribution with minimum $54$, mode $55$ and maximum $56$.
A $10^5$ MC sample  gives an estimation of the failure probability $\hat{P}=8.6\times10^{-4}$.

\subsubsection{FORM}
The algorithm FORM converges to a design point $\left(1.72,-2.70,0.55,-0.18\right)$ in $52$ function calls, giving an approximate probability of $\hat{P}_{FORM}=5.8\times10^{-4}$. The importance factors are displayed in Table \ref{tab:IFCrue}. 

{\footnotesize }
\begin{table}[h]
\begin{centering}
{\footnotesize }%
\begin{tabular}{|c|c|c|c|c|}
\hline 
{\footnotesize Variable} & {\footnotesize $Q$} & {\footnotesize $K_{s}$ } & {\footnotesize $Z_{v}$ } & {\footnotesize $Z_{m}$ }\tabularnewline
\hline 
\hline 
{\footnotesize Importance factor} & {\footnotesize $0.246$} & {\footnotesize $0.725$} & {\footnotesize $0.026$} & {\footnotesize $0.003$}\tabularnewline
\hline 
\end{tabular}
\par\end{centering}{\footnotesize \par}

{\footnotesize \caption{\label{tab:IFCrue}{\footnotesize Importance factors for flood case}}
}
\end{table}
{\footnotesize \par}

FORM assesses that $K_{s}$ is of extremely high influence, followed by $Q$ that is of medium influence. $Z_{v}$ has a very weak influence and $Z_{m}$ is negligible. It can be noticed that the estimated failure probability is twice as small as the one estimated with crude MC, but remains in the same order of magnitude.

\subsubsection{Sobol indices}
The first-order and total indices are displayed in Table \ref{tab:SobolCrue}. 
It can be seen that the estimates of some indices are negative despite the fact that Sobol's indices are theoretically positive. The estimation can indeed produce negative results for values close to 0.

{\footnotesize }
\begin{table}[h]
\begin{centering}
{\footnotesize }%
\begin{tabular}{|c|c|c|c|c|c|c|}
\hline 
{\footnotesize Sobol' Index} & {\footnotesize $S_{1}$} & {\footnotesize $S_{2}$} & {\footnotesize $S_{3}$} & {\footnotesize $S_{T1}$} & {\footnotesize $S_{T2}$} & {\footnotesize $S_{T3}$}\tabularnewline
\hline 
\hline 
{\footnotesize Mean} & {\footnotesize $0.0234$} & {\footnotesize $0.0099$} & {\footnotesize $0.0667$} & {\footnotesize $0.8158$} & {\footnotesize $0.6758$} & {\footnotesize $0.9299$}\tabularnewline
\hline 
{\footnotesize C.o.v.} & {\footnotesize $0.0072$} & {\footnotesize $0.0051$} & {\footnotesize $0.0095$} & {\footnotesize $0.0156$} & {\footnotesize $0.0216$} & {\footnotesize $0.0094$}\tabularnewline
\hline 
\end{tabular}
\par\end{centering}{\footnotesize \par}

{\footnotesize \caption{\label{tab:SobolIshi}{\footnotesize Sobol' indices for Ishigami function}}
}
\end{table}
{\footnotesize \par}

Considering the first order indices, $Z_{v}$ and $Z_{m}$ are of null influence on their own. $Q$ is considered to have a minimal influence ($1\%$ of the variance of the indicator function) by itself, and $K_{s}$ explains $24\%$ of the variance on its own. When considering the total indices, it can be noticed that both $Z_{v}$ and $Z_{m}$ have a weak impact on the failure probability. On the other hand, $Q$ has a major influence on the failure probability. $K_{s}$ total index is close to one, therefore $K_{s}$ explains (with or without any interaction with other variables) almost all the variance of the failure function.

Let us compare the informations provided by the Sobol' indices with the information provided by the importance factors. One cannot conclude from the total Sobol' indices that $Z_m$ is not influent whereas the importance factors assess that this variable is of negligible influence. Additionally, the total Sobol' index associated to $K_s$ and $Q$ state that both these variables are of high influence whereas the importance factors state that $K_s$ is of high influence and $Q$ is of medium influence. 

\subsubsection{Density modification based reliability indices }\label{IBSMflood}
The method presented throughout this article is applied on the flood case. Only the mean shifting will be applied here. The modified pdf are given in Table \ref{tab:TableRecap} (Appendix \ref{TableE})and a numerical trick is used to deal with truncated distributions, as stressed in Appendix \ref{Numtricks}. One can notice that the different inputs follow various distributions (unlike the other examples), thus the question of "equivalent" perturbation arises. It will be discussed further in Section \ref{sec:Discussion}. Here the choice has been made to shift the mean relatively to the standard deviation, hence including the spread of the various inputs in their respective perturbation. So for any input, the original distribution is perturbed so that its mean is the original's one plus $\delta$ times its standard deviation, $\delta$ ranging from $-1$ to $1$ with $40$ points.

\begin{figure}
\includegraphics[width=1\textwidth]{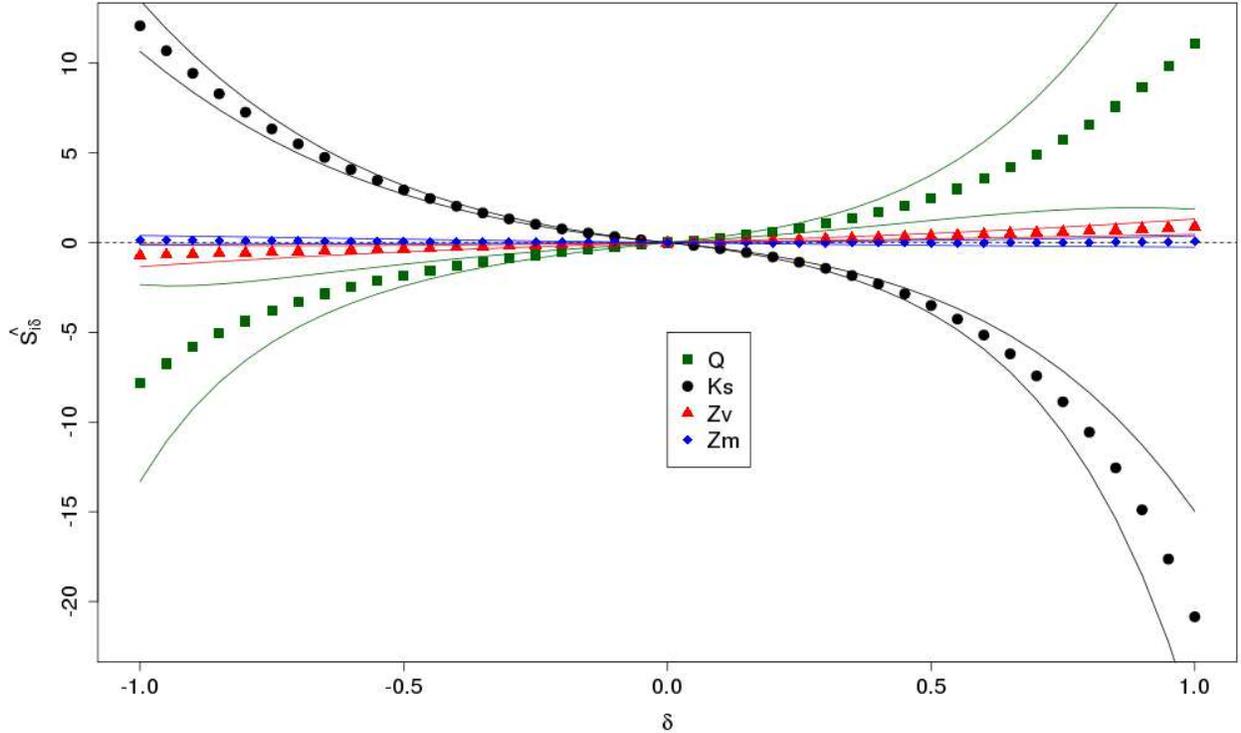}
\caption{\label{fig:IBSMDCrue}
Estimated indices $\widehat{S_{i\delta}}$ for the flood case  with a mean perturbation}
\end{figure}

Figure \ref{fig:IBSMDCrue} assesses that an increasing of the mean of the inputs increases the failure probability slightly for  $Z_{v}$, strongly for  $Q$, and diminishes it slightly for  $Z_{m}$ and strongly for  $K_{s}$. This goes the opposite way when decreasing the mean. In terms of absolute modification, $K_{s}$ and $Q$ are of same magnitude, even if  $K_{s}$ has a slightly stronger impact. On the other hand, the effects of mean perturbation on $Z_{m}$ and $Z_{v}$ are negligible. 
The CI associated to $Q$ and $K_s$ are well separated from the others, except in a $\delta=-.3$ to $.3$ zone. The CI associated to $Z_v$ and $Z_m$ overlap. Thus even though the indices seem to have different value, it is not possible to conclude with certainty about the influence of those variables.

%%%%%%%%%%%%%%%%%%%%%%%%%%%%%%%%%%%%%%%%%%%%%%%%%%%%%%%%%%%%%%%%%%%%%%%%%%%%%%%%%%%%%%%%%%%%%%%%%%%%%%%%%%%%%%%
  
\section{Discussion \label{sec:Discussion}}
 
\subsection{Conclusion on the method}
The method presented in this paper gives relevant complementary information in addition of traditional SA methods applied to a reliability problem. Traditional SA methods provide variable ranking, whereas the proposed method provides an indication on the variation in the probability of failure given the variation of  parameter $\delta$. This is useful when the practitioner is interested on which configurations of the problem lead to an increase of the failure probability. This might also be used to assess the conservatism of a problem, if  every variations of the input lead to decrease in the probability of failure.  Additionally, it  has two advantages: 
\begin{itemize}
\item The ability for the user to set the most adapted constraints considering his/her problem/objective,
\item The MC framework allowing to use previously done function calls, thus limiting the CPU cost of the SA, and allowing the user to test several perturbations.
\end{itemize} 
 We argue that with an adapted perturbation, this method can fulfil the three presented engineer’s objective.

\subsection{Equivalent perturbation}
The question of "equivalent" perturbation arises from cases where all inputs are not identically distributed. Indeed, problems may emerge when some inputs are defined on infinite intervals and when other inputs are defined on finite intervals (such as uniform distributions). For instance, consider a two-dimensional model with one Gaussian distribution and one uniform distribution as inputs. Thus, a mean shift will be a translation for the first input, whereas it will lead to a Dirac distribution in one endpoint for the other input. Hence, a mean shift cannot be considered as an "equivalent" perturbation. 
A "relative mean shift" seems promising idea. But if we consider a model with two Gaussian inputs of equal variance 1 and of mean respectively 1 and 10000. Then, a relative mean shift of 10\% will result in Gaussian distributions with mean respectively 1.1 and 11000, and still variance 1. This counter-example shows that relative mean shift might not be an adequate perturbation in terms of "equivalence".

\subsection{Support pertubation}
In the examples given throughout this paper, the perturbations of the inputs left the support of those variables unperturbed. However, the practitioner might be interested  in the sensitivity on the boundaries of the support. The proposed method will be applied with support perturbations in further tests. However, we stress that given the estimation method (reverse importance sampling), it is mandataory that the support of the perturbed density is included in the support of the original density. Thus one cannot perturb the inputs so that the perturbed support is wider than the original one. 

\subsection{Further work}
Finally, these first results provide some avenues for future research: 
\begin{itemize}
\item Adapting the estimator of the indices $S_{i\delta}$ in term of variance reduction and of number of function calls. Further work will be made with importance sampling methods, and possibly subset methods. The use of sequential methods \cite{bect2011sequential} may also be tested.  
\item Second, there is a need to find a way to perturb "equivalently" several distributions of different natures. To answer this problem, perturbations that are not based upon moment constraints can be proposed.

For instance, quantile constraints might be considered. As far as we noticed, in most cases the values of the input leading to the failure event comes from the tails of the input distributions. What if these tails were badly modeled ? Therefore a perturbation based on the quantiles is proposed. More precisely, suppose that the weigth of the left tail is increased. That is to say that the value $q_{5\%}$ becomes for the modified density, for instance the $\delta$ quantile. This writes:
$$\int1_{]-\infty;q_{5\%}]}(x)f_{\mbox{\tiny mod}}(x)dx=\delta.$$
We believe that this kind of perturbation is equivalent with respect to inputs of different natures.

A perturbation based on an entropy constraint might also be proposed. The differential entropy of a distribution can be seen as a quantification of uncertainty \cite{auder2009global}. Thus an example of (non-linear) constraint on the entropy can be:
$$
-\int f_{X_{i\delta}}(x)\log f_{X_{i\delta}}(x)dx=-\delta\int f_{X_{i}}(x)\log f_{X_{i}}(x)dx.
$$
Yet further computations have to be made to obtain a tractable solution of the KL minimization problem under the above constraint. 
\end{itemize}

%%%%%%%%%%%%%%%%%%%%%%%%%%%%%%%%%%%%%%%%%%%%%%%%%%%%%%%%%%%%%%%%%%%%%%%%%%%%%%%%%%%%%%%%%%%%%%%%%%%%%%%%%%%%%%%
\section*{Acknoweldgements \label{sec:ack}}

Part of this work has been backed by French National Research Agency (ANR)
through COSINUS program (project COSTA BRAVA noANR-09-COSI-015).
We thank Dr. Daniel Busby (IFP EN) for several discussions and Emmanuel Remy (EDF R\&D) for proofreading.
We also thank an anonymous reviewer and the associated editor for their helpful comments.
All the statistical parts of this work have been performed within the R environment, including the \emph{mistral} and \emph{sensitivity} packages. 

%%%%%%%%%%%%%%%%%%%%%%%%%%%%%%%%%%%%%%%%%%%%%%%%%%%%%%%%%%%%%%%%%%%%%%%%%%%%%%%%%%%%%%%%%%%%%%%%%%%%%%%%%%%%%%%
%\section*{References \label{sec:ref}}
%
\bibliographystyle{plain}
\bibliography{Article}

%%%%%%%%%%%%%%%%%%%%%%%%%%%%%%%%%%%%%%%%%%%%%%%%%%%%%%%%%%%%%%%%%%%%%%%%%%%%%%%%%%%%%%%%%%%%%%%%%%%%%%%%%%%%%%%
\section*{Appendices}
\appendix

%%%%%%%%%%%%%%%%%%%%%%%%%%%%%%%%%%%%%%%%%%%%%%%%%%%%%%%%%%%%%%%%%%%%%%%%%%%%%%%%%%%%%%%%%%%%%%%%%%%%%%%%%%%%%%%
\section{Computation of Lagrange multipliers}\label{appLag}% $\boldsymbol{\lambda}$}
%%%%%%%%%%%%%%%%%%%%%%%%%%%%%%%%%%%%%%%%%%%%%%%%%%%%%%%%%%%%%%%%%%%%%%%%%%%%%%%%%%%%%%%%%%%%%%%%%%%%%%%%%%%%%%%

%We hereby note $\boldsymbol{\lambda}$ the optimal Lagrange vector
%-which is wanted- and $\mathbf{l}=(l_{1},...,l_{K})$ a vector of
%constants.
Let $H$ be the Lagrange function: 
\[
H(\boldsymbol{\lambda})=\psi_i(\boldsymbol{\lambda})-\sum_{k=1}^{K}\lambda_{k}\delta_{k}
.\]
Thus, using the results of \cite{csiszar1975divergence}, one has
\[
\boldsymbol{\lambda}^*=\arg\min H(\boldsymbol{\lambda})
.\]
The expression of the gradient of $H$ with respect to the $j^{\mbox{th}}$ variable is
\[
\nabla_{j}H(\boldsymbol{\lambda})=\frac{\int g_{j}(x)f_{i}(x)\exp(\sum_{k=1}^{K}\lambda_{k}g_{k}(x))dx}{\exp\psi_i(\boldsymbol{\lambda})}-\delta_{j}.
\]
Similarly, the expression of the second derivative  of $H$ with respect to
the $h^{\mbox{th}}$ and the $j^{\mbox{th}}$ variables is
\begin{align*}
&D_{hj}H(\boldsymbol{\lambda}) = \frac{\int g_{h}(x)g_{j}(x)f_{i}(x)\exp(\sum_{k=1}^{K}\lambda_{k}g_{k}(x))dx}{\exp\psi_i(\boldsymbol{\lambda})}\\
 &-\frac{\int g_{j}(x)f_{i}(x)\exp(\sum_{k=1}^{K}\lambda_{k}g_{k}(x))dx}{\exp\psi_i(\boldsymbol{\lambda})}\frac{\int g_{h}(x)f_{i}(x)\exp(\sum_{k=1}^{K}\lambda_{k}g_{k}(x))dx}{\exp\psi_i(\boldsymbol{\lambda})}.
\end{align*}
This method has been used in this paper for computing the optimal vector $\boldsymbol{\lambda}^*$ when a variance shifting was applied.
 The integrals were evaluated with Simpson's rule.

%%%%%%%%%%%%%%%%%%%%%%%%%%%%%%%%%%%%%%%%%%%%%%%%%%%%%%%%%%%%%%%%%%%%%%%%%%%%%%%%%%%%%%%%%%%%%%%%%%%%%%%%%%%%%%%
\section{Proofs of the NEF properties}{\label{appedix:proofsofnef}}
%%%%%%%%%%%%%%%%%%%%%%%%%%%%%%%%%%%%%%%%%%%%%%%%%%%%%%%%%%%%%%%%%%%%%%%%%%%%%%%%%%%%%%%%%%%%%%%%%%%%%%%%%%%%%%%
In this Appendix, the details of the calculus for the Proposition \ref{prop:NEF} are provided. 
\paragraph*{NEF specificities : }
If the original density $f_{i}(x)$ is a NEF, then under a set of
$K$ linear constraints on $f(x)$, one has :
\[
f(x)=b(x)\exp\left[x\theta-\eta(\theta)\right],
\]
thus :
\[
f_{\delta}(x)=f(x)\exp\left[\sum_{k=1}^{K}\lambda_{k}g_{k}(x)-\psi(\mathbf{\boldsymbol{\lambda}})\right]
\]
The regularization constant from (\ref{eq:phi}) can be written as: 
\begin{equation}
\psi(\boldsymbol{\lambda}) = \log\int b(x)\exp\left[x\theta+\sum_{k=1}^{K}\lambda_{k}g_{k}(x)-\eta(\theta)\right]dx\label{eq:reg}
\end{equation}
If the integral on (\ref{eq:reg}) is finite, $f_{\delta}$ exists and
is a density.
\paragraph*{Mean shifting}
With a  single constraint formulated as in (\ref{eq:meanshift}), (\ref{eq:reg}) becames :
\begin{eqnarray*}
\psi(\boldsymbol{\lambda}) & = &\log\int b(x)\exp\left[x\theta+\lambda x-\eta(\theta)\right]dx \\
& = &\log\int b(x)\exp\left[x\left(\theta+\lambda\right)
-\eta(\theta)+\eta(\theta+\lambda)-\eta(\theta+\lambda)\right]dx 
\end{eqnarray*}
if $\eta(\theta+\lambda)$ is well defined.
\begin{eqnarray*}
\psi(\boldsymbol{\lambda}) & = & \left(\eta(\theta+\lambda)-\eta(\theta)\right)+\log\left[\int b(x)\exp\left[x\left(\theta+\lambda\right)-\eta(\theta+\lambda)\right]\right]dx\\
& = & \eta(\theta+\lambda)-\phi(\theta)
\end{eqnarray*}
since 
\[
b(x)\exp\left[x\left(\theta+\lambda\right)-\eta(\theta+\lambda)\right]=f_{\theta+\lambda}(x)
\]
with notation from (\ref{eq:notationNEF}), is a density of integral 1. Thus
\begin{eqnarray*}
f_{\delta}(x) & = & b(x)\exp\left[x\theta-\phi(\theta)\right]\exp\left[\lambda x-\eta(\theta+\lambda)+\eta(\theta)\right] \\
& = & b(x)\exp\left[x\left[\theta+\lambda\right]- \eta(\theta+\lambda)\right]=f_{\theta+\lambda}(x)
\end{eqnarray*}
Thus the mean shifting of a NEF of CDF $\eta(.)$ results in another
NEF with mean $\eta'(\theta+\lambda)=\delta$ (constraint) and variance
$\eta''(\theta+\lambda)$.
\paragraph*{Variance shifting}
With a  single constraint formulated as in (\ref{eq:cstrvariance}), using (\ref{eq:reg}), the new distribution has for density:
\[
f_{\delta}(x)=b(x)\exp\left[x\theta+x\lambda_{1}+x^{2}\lambda_{2}-\psi(\mathbf{\boldsymbol{\lambda}})-\eta(\theta)\right]
\]
Since $\boldsymbol{\lambda}$ is known or computed, and $\theta$
is also known, consider the variable change $
z=\sqrt{\lambda_{2}}x$ assuming $\lambda_{2}$ is strictly positive (the variable change is $z=\sqrt{-\lambda_{2}}x$ if $\lambda_{2}$ is strictly negative).
Thus, 
\begin{eqnarray*}
f_{\delta}(x) & = & b(\frac{z}{\sqrt{\lambda_{2}}})\exp\left[z^{2}\right]\exp\left[\frac{z}{\sqrt{\lambda_{2}}}\left(\theta+\lambda_{1}\right)-\psi(\mathbf{\boldsymbol{\lambda}})-\eta(\theta)\right] \\
%& = & c(z)\exp\left[z\frac{\left(\theta+\lambda_{1}\right)}{\sqrt{\lambda_{2}}}-\eta(\theta)+\eta\left(\frac{\left(\theta+\lambda_{1}\right)}{\sqrt{\lambda_{2}}}\right)-\eta\left(\frac{\left(\theta+\lambda_{1}\right)}{\sqrt{\lambda_{2}}}\right)\right]\frac{1}{\exp\left[\psi(\mathbf{\boldsymbol{\lambda}})\right]} \\
& = & \exp\left[\eta\left(\frac{\left(\theta+\lambda_{1}\right)}{\sqrt{\lambda_{2}}}\right)-\eta(\theta)-\psi(\mathbf{\boldsymbol{\lambda}})\right]c(z)\exp\left[z\frac{\left(\theta+\lambda_{1}\right)}{\sqrt{\lambda_{2}}}-\eta\left(\frac{\left(\theta+\lambda_{1}\right)}{\sqrt{\lambda_{2}}}\right)\right]\label{eq:fin1}
\end{eqnarray*}
with 
\[
c(z)=b(\frac{z}{\sqrt{\lambda_{2}}})\exp\left[z^{2}\right].
\]
By (\ref{eq:phi}),
\begin{eqnarray*}
\psi(\boldsymbol{\lambda}) & = & \log\int b(x)\exp\left[x\theta+x\lambda_{1}+x^{2}\lambda_{2}-\eta(\theta)\right]dx \\
& = &\log\int b(\frac{z}{\sqrt{\lambda_{2}}})\exp\left[z^{2}\right]\exp\left[\frac{\left(\theta+\lambda_{1}\right)}{\sqrt{\lambda_{2}}}z-\eta(\theta)+\eta\left(\frac{\left(\theta+\lambda_{1}\right)}{\sqrt{\lambda_{2}}}\right)-\eta\left(\frac{\left(\theta+\lambda_{1}\right)}{\sqrt{\lambda_{2}}}\right)\right]dx \\
& = & \left(\eta\left(\frac{\left(\theta+\lambda_{1}\right)}{\sqrt{\lambda_{2}}}\right)-\eta(\theta)\right)+\log\int c(z)\exp\left[\frac{\left(\theta+\lambda_{1}\right)}{\sqrt{\lambda_{2}}}z-\eta\left(\frac{\left(\theta+\lambda_{1}\right)}{\sqrt{\lambda_{2}}}\right)\right]dx\\
& = & \eta\left(\frac{\left(\theta+\lambda_{1}\right)}{\sqrt{\lambda_{2}}}\right)-\eta(\theta)\label{eq:fin2}
\end{eqnarray*}

Thus one has : 
\[
f_{\delta}(x)=c(z)\exp\left[z\frac{\left(\theta+\lambda_{1}\right)}{\sqrt{\lambda_{2}}}-\eta\left(\frac{\left(\theta+\lambda_{1}\right)}{\sqrt{\lambda_{2}}}\right)\right]
\]
thus the variance shifting of a NEF results in another NEF parameterized by $\frac{\left(\theta+\lambda_{1}\right)}{\sqrt{\lambda_{2}}}$.

%%%%%%%%%%%%%%%%%%%%%%%%%%%%%%%%%%%%%%%%%%%%%%%%%%%%%%%%%%%%%%%%%%%%%%%%%%%%%%%%%%%%%%%%%%%%%%%%%%%%%%%%%%%%%%%
\section{Proofs of asymptotic properties}\label{proofs}
%%%%%%%%%%%%%%%%%%%%%%%%%%%%%%%%%%%%%%%%%%%%%%%%%%%%%%%%%%%%%%%%%%%%%%%%%%%%%%%%%%%%%%%%%%%%%%%%%%%%%%%%%%%%%%%

\subsection{Proof of Lemma \ref{lemma1}} \label{proofoflemma1}
Under assumption {\bf(i)}, we have
$$
%& = & \int_{\mbox{Supp}(f_{i})} f_{i\delta}(x) \ d\mathbf{x}, \\
\int_{\text{\scriptsize Supp}(f_{i\delta})} \1_{\{G(\mathbf{\mathbf{x}})<0\}}\frac{f_{i\delta}(x_i)}{f_i(x_i)}f(\mathbf{x}) \ d\mathbf{x}  \leq  \int_{\text{\scriptsize Supp}(f_{i\delta})} f_{i\delta}(x_i) \ dx_i \ = \ 1.
$$

So that, the strong LLN  may be applied to $\hat{P}_{i\delta N}$. Defining  
\begin{equation}\label{eq:expression_sigma}
\sigma^{2}_{i\delta} = \mbox{Var}\left[\1_{\{G(\mathbf{X})<0\}}\frac{f_{i\delta}(X_{i})}{f_{i}(X_{i})}\right],
\end{equation}
 one has
\begin{eqnarray*}
{\sigma}^{2}_{i\delta} & = & \int_{\text{\scriptsize Supp}(f_{i})} \1_{\{G(\mathbf{x})<0\}}  \frac{f_{i\delta}^2(x_{i})}{f_{i}(x_{i})} \prod\limits_{j\neq i}^d f_{j}(x_{j}) \ d{\bf x} -P_{i\delta}^2\ < \ \infty \ \ \ \ \text{under Condition {\bf(ii)}.}
\end{eqnarray*}
Therefore the CLT applies: 
\begin{eqnarray*}
\sqrt{N}{\sigma}^{-1}_{i\delta}\left(\hat{P}_{i\delta N}- P_{i\delta}\right)  \xrightarrow[]{{\cal{L}}}  {\cal{N}}(0,1)\;.
\end{eqnarray*}
%With $\hat{\sigma}^{2}_{i\delta N}$ the standard Monte Carlo estimator of ${\sigma}^{2}_{i\delta}$, and $\E[\hat{\sigma}^{2}_{i\delta N}]<\infty$ 

Under assumption {\bf(ii)}, the strong LLN applies to $\hat{\sigma}^{2}_{i\delta N}$.  So that, the final result is straightforward using Slutsky's lemma.

\subsection{Proof of Proposition \ref{CTL.correlated}\label{A2}} First, note that 
\begin{eqnarray*}
%\text{Cov}(\hat{P}_{N},\hat{P}_{i\delta N}) & =& \mathbb{E}\left[\widehat{P}\widehat{P_{i\delta}}\right]-\mathbb{E}\left[\hat{P}_{N}\right]\mathbb{E}\left[\hat{P}_{i\delta N}\right], \\ 
%& =& 
\mathbb{E}\left[\widehat{P}\widehat{P_{i\delta}}\right]-PP_{i\delta} 
& =  & \mathbb{E}\left[\frac{1}{N^{2}}\left(\sum_{n=1}^{N}\1_{\{G(\mathbf{x}^{n})<0\}}\right)\left(\sum_{n=1}^{N}{\1}_{\{G(\mathbf{x}^{n})<0\}}\frac{f_{i\delta}(x_{i}^{n})}{f_{i}(x_{i}^{n})}\right)\right]-PP_{i\delta} \\
& = & \frac{1}{N^{2}}\mathbb{E}\left[\sum_{n=1}^{N}\left[\1_{\{G(\mathbf{x}^{n})<0\}}\right]^{2}\frac{f_{i\delta}(x_{i}^{n})}{f_{i}(x_{i}^{n})}+\sum_{n=1}^{N}\sum_{j\neq i}^{N}{\1}_{\{G(\mathbf{x}^{n})<0\}}\1_{\{G(\mathbf{x}^{j})<0\}}\frac{f_{i\delta}(x_{i}^{j})}{f_{i}(x_{i}^{j})}\right]\\
 && -PP_{i\delta} \\
 & = & \frac{1}{N^{2}}\left[NP_{i\delta}+N\left(N-1\right)PP_{i\delta}\right]-PP_{i\delta} \\
& = & \frac{1}{N}\left(P_{i\delta}-PP_{i\delta}\right)\;.
\end{eqnarray*}
Assuming the conditions under which Lemma 1 is true, the bivariate CLT follows with
\begin{eqnarray*}\label{eq:expression_Sigma}
\Sigma_{i\delta} & = & \left(\begin{array}{cc}
P(1-P) & P_{i\delta}(1-P)\\
P_{i\delta}(1-P) & \sigma_{i\delta}^2\
\end{array}\right)\;.
\end{eqnarray*}
Each term of this matrix can be consistently estimated, using the results in Lemma 1 and Slutsky's lemma.

%%%%%%%%%%%%%%%%%%%%%%%%%%%%%%%%%%%%%%%%%%%%%%%%%%%%%%%%%%%%%%%%%%%%%%%%%%%%%%%%%%%%%%%%%%%%%%%%%%%%%%%%%%%%%%%
\section{Summary Table with modified distributions for mean shift}\label{TableE}
%%%%%%%%%%%%%%%%%%%%%%%%%%%%%%%%%%%%%%%%%%%%%%%%%%%%%%%%%%%%%%%%%%%%%%%%%%%%%%%%%%%%%%%%%%%%%%%%%%%%%%%%%%%%%%%

\begin{sidewaystable}
    
	\centering
	\footnotesize
	\caption{Modified distributions for mean shifting.Note that $\Phi(.)$is the cdf of the standard normal distribution, and $\phi(.)$ is its pdf.}

	{\begin{tabular}{cccc}
	\hline 
	{Original distribution} & {Modified distribution} & {Modified pdf $f_{i\delta}$} & {Link between $\lambda^{*}$ and $\delta$}
	\tabularnewline
	
	\hline 
	\hline 
	
	{ $\mbox{{NEF}}(\theta)$} &
	{ $NEF(\theta+\lambda^{*})$} &
	{ $f_{i\delta}(x_{i})=b(x_{i})\exp\left[x_{i}\left[\theta+\lambda^{*}\right]-\eta(\theta+\lambda^{*})\right]$} & 
	{ $\eta'(\theta+\lambda^{*})=\delta$}
	\tabularnewline
	
	{ Special case of NEF: $\mathcal{N}(\mu,\sigma)$} & 
	{ $\mathcal{N}(\delta,\sigma)$} & 
	{ $f_{i\delta}(x_{i})=\frac{1}{\sigma\sqrt{2\pi}}\exp\left[-\frac{1}{2}\left(\frac{x_{i}-\delta}{\sigma}\right)^{2}\right]$} &
	{ $\lambda^{*}=\frac{\delta-\mu}{\sigma^{2}}$}
	\tabularnewline
	
	{ Uniform distribution: $\mathcal{U}_{[a,b]}$} & { $\propto$ truncated exponential } & 
	{ $ f_{i\delta}(x_{i})=\frac{\lambda^{*}}{e^{\lambda^{*}b}-e^{\lambda^{*}a}}1_{[a,b]}(x_{i})e^{\lambda x_{i}}$} & 
	{ $\delta=\frac{1}{\left(b-a\right)}\frac{e^{\lambda^{*}b}\left(\lambda^{*}b-1\right)+e^{\lambda^{*}a}\left(1-\lambda^{*}a\right)}{\lambda^{*}\left(e^{\lambda^{*}b}-e^{\lambda^{*}a}\right)}$}
	\tabularnewline
	
	{ Left Tr Gaussian $\mathcal{N}_{T}(\mu,\sigma,a)$} &
	{ $\mathcal{N}_{T}(\mu+\sigma^{2}\lambda^{*},\sigma,a)$} &
	{ $f_{i\delta}(x_{i})=\frac{1_{[a,+\infty[}(x_{i})}{1-F(a)}\frac{1}{\sigma\sqrt{2\pi}}\exp\left[-\frac{1}{2}\left(\frac{x_{i}-\mu-\sigma^{2}\lambda^{*}}{\sigma}\right)^{2}\right]$} &
	{ $\delta=\mu+\sigma^{2}\lambda^{*}-\sigma\frac{\phi\left(\frac{a-(\mu+\sigma^{2}\lambda^{*})}{\sigma}\right)}{1-\Phi\left(\frac{a-(\mu+\sigma^{2}\lambda^{*})}{\sigma}\right)}$}
	\tabularnewline
	
	{ Triangle $\mathcal{T}(a,b,c)$} &
	{ --} &
	{ $f_{i\delta}(x_{i})=\exp(x_{i}\lambda^{*}-\psi(\lambda^{*}))f(x_{i})$} &
	{ $\delta=\frac{(a-\frac{1}{\lambda^{*}})e^{\lambda^{*}a}(b-c)+(b-\frac{1}{\lambda^{*}})e^{\lambda^{*}b}(c-a)+(c-\frac{1}{\lambda^{*}})e^{\lambda^{*}c}(a-b)}{e^{\lambda^{*}a}(b-c)+e^{\lambda^{*}b}(c-a)+e^{\lambda^{*}c}(a-b)}$}
	\tabularnewline
	
	{ Left Tr Gumbel $\mathcal{G}_{T}(\mu,\beta,a)$} &
	{ --} &
	{ $f_{i\delta}(x_{i})=\exp(x_{i}\lambda^{*}-\psi(\lambda^{*}))f(x_{i})$} &
	{ $	\begin{array}{c}
			\delta=\frac{M_{Y}'(\lambda^{*})-\int_{-\infty}^{a}yf_{Y}(y)\exp\left[\lambda^{*}y\right]dy}{M_{Y}(\lambda^{*})-\int_{-\infty}^{a}f_{Y}(y)\exp\left[\lambda^{*}y\right]dy}\mbox{ with :}\\
			M_{Y}(\lambda^{*})=\Gamma\left(1-\beta\right)\exp\left[\lambda^{*}\mu\right]\\
			M_{Y}'(\lambda^{*})=\Gamma\left(1-\beta\right)\exp\left[\lambda^{*}\mu\right]\left[\mu-\beta\digamma^{(0)}(1-\lambda^{*})\right]
		\end{array}$}
	\tabularnewline
	\hline
	 
\end{tabular}}
\label{tab:TableRecap}
\end{sidewaystable}

%%%%%%%%%%%%%%%%%%%%%%%%%%%%%%%%%%%%%%%%%%%%%%%%%%%%%%%%%%%%%%%%%%%%%%%%%%%%%%%%%%%%%%%%%%%%%%%%%%%%%%%%%%%%%%%
\section{Numerical trick to work with truncated distribution}\label{Numtricks}
%%%%%%%%%%%%%%%%%%%%%%%%%%%%%%%%%%%%%%%%%%%%%%%%%%%%%%%%%%%%%%%%%%%%%%%%%%%%%%%%%%%%%%%%%%%%%%%%%%%%%%%%%%%%%%%

In the case where a mean shifting is considered on a left truncated
distribution. We present a tip that can help to compute $\boldsymbol{\lambda}^{*}$.

The studied trucated variable $Y_{T}$ has distribution $f_{YT}.$
Let us denote $Y\sim f_{Y}$ the corresponding non-truncated distribution.
The truncation occurs for some real value $a$. This truncation may
happen for some physical modelling reason. One has:
\[
f_{YT}(y)=\frac{1}{1-F(a)}1_{[a,+\infty[}(y)f_{Y}(y).
\]
The formal definition of $M_{YT}(\boldsymbol{\lambda})$ the mgf of
$Y_{T}$ for some $\boldsymbol{\lambda}$ is:
\[
M_{YT}(\boldsymbol{\lambda})=\frac{1}{1-F_{Y}(a)}\int_{a}^{+\infty}f_{Y}(y)\exp\left[\boldsymbol{\lambda}y\right]dy.
\]
Let us recall that we are looking for $\boldsymbol{\lambda}^{*}$such
as:
\begin{equation}
\delta=\frac{M_{YT}'(\boldsymbol{\lambda}^{*})}{M_{YT}(\boldsymbol{\lambda}^{*})}=\frac{\int_{a}^{+\infty}yf_{Y}(y)\exp\left[\boldsymbol{\lambda}y\right]dy}{\int_{a}^{+\infty}f_{Y}(y)\exp\left[\boldsymbol{\lambda}y\right]dy}.\label{eq:deltaexpre1}
\end{equation}

When the expression does not take a practical form, one can use numerical
integration to estimate the integral terms. Unfortunately, for some
heavy tailed distribution (for instance Gumbel distribution), this
numerical integration might be complex or not possible. This is due
to the multiplication by an exponential of $y$. The following tip
helps to avoid such problems. Denoting $M_{Y}(\boldsymbol{\lambda})$
the mgf of the non-truncated distribution, one can remark
that:
\[
M_{Y}(\boldsymbol{\lambda})=\int_{-\infty}^{+\infty}f_{Y}(y)\exp\left[\boldsymbol{\lambda}y\right]dy=\int_{-\infty}^{a}f_{Y}(y)\exp\left[\boldsymbol{\lambda}y\right]dy+\int_{a}^{+\infty}f_{Y}(y)\exp\left[\boldsymbol{\lambda}y\right]dy
\]
Thus another expression for $M_{YT}(\boldsymbol{\lambda})$ is:
\[
M_{YT}(\boldsymbol{\lambda})=\frac{1}{1-F_{Y}(a)}\left[M_{Y}(\boldsymbol{\lambda})-\int_{-\infty}^{a}f_{Y}(y)\exp\left[\boldsymbol{\lambda}y\right]dy\right].
\]
The integral term is much smaller in the left heavy tailed distribution
case. Therefore the numerical integration (for instance using Simpson's
method) is much more precise or became possible. 

The same goes for $M_{YT}'(\boldsymbol{\lambda})$ which has alternative
expression:
\[
M_{YT}'(\boldsymbol{\lambda})=\frac{1}{1-F_{Y}(a)}\left[M_{Y}'(\boldsymbol{\lambda})-\int_{-\infty}^{a}yf_{Y}(y)\exp\left[\boldsymbol{\lambda}y\right]dy\right].
\]
Finally, another form of \ref{eq:deltaexpre1} is:
\begin{equation}
\delta=\frac{M_{Y}'(\boldsymbol{\lambda})-\int_{-\infty}^{a}yf_{Y}(y)\exp\left[\boldsymbol{\lambda}y\right]dy}{M_{Y}(\boldsymbol{\lambda})-\int_{-\infty}^{a}f_{Y}(y)\exp\left[\boldsymbol{\lambda}y\right]dy}.\label{eq:deltaexpre1-1}
\end{equation}

This alternative expression may lead to more precise estimations of
$\boldsymbol{\lambda}^{*}$ when $M_{Y}(\boldsymbol{\lambda})$ and
$M_{Y}'(\boldsymbol{\lambda})$ are known (which is the case for most
usual distribution) since the integral term are much smaller than
in the first expression . A reference to this Appendix is made in the summary table \ref{tab:TableRecap}.

\end{document}